\documentclass{amsart}
\usepackage{latexsym}
\usepackage{amssymb}
\usepackage{amsmath}
\usepackage{amsthm}
\usepackage{hyperref}
\usepackage{xcolor}

\newcommand{\CC}{{\mathbb C}}

\newcommand{\RR}{{\mathbb R}}

\newcommand{\SSS}{{\mathbb S}}

\newcommand{\grad}{\mathrm{grad}}

\newcommand{\scal}{\mathrm{Scal}}

\newcommand{\tr}{\mathrm{tr}}

\newcommand{\lgra}{\longrightarrow}

\newtheorem{example}{Examples}[section]
\newtheorem{thm}{Theorem}[section]
\newtheorem{lemma}[thm]{Lemma}
\newtheorem{prop}[thm]{Proposition}
\newtheorem{cor}[thm]{Corollary}
\newtheorem{remark}[thm]{Remark}
\newtheorem{remarks}[thm]{Remarks}
\newtheorem{definition}[thm]{Definition}
\newtheorem{notation}[thm]{Notation}
\newtheorem{exabout:ample}[thm]{Example}

\begin{document}
\title{$f$-Biharmonic and bi-$f$-harmonic submanifolds of generalized space forms}
\author[J. ROTH]{Julien Roth}
\address[J. ROTH]{Laboratoire d'Analyse et de Math\'ematiques Appliqu\'ees, UPEM-UPEC, CNRS, F-77454 Marne-la-Vall\'ee}
\email{julien.roth@u-pem.fr}
\author[A. UPADHYAY]{Abhitosh Upadhyay}
 \address[A. UPADHYAY]{Harish Chandra Research Institute, Chhatnag Road, Jhunsi, Allahabad, India, 211019}
\email{abhi.basti.ipu@gmail.com, abhitoshupadhyay@hri.res.in}

\maketitle

\begin{abstract}
We study $f$-biharmonic and bi-$f$-harmonic submanifolds in both generalized \linebreak complex and Sasakian space forms. We prove necessary and sufficient condition for \linebreak $f$-biharmonicity and bi-$f$-harmonicity in the general case and many particular cases. Some non-existence results are also obtained.
\end{abstract}

\section{\textbf{Introduction}}
Harmonic maps between two Riemannian manifolds $(M^{m},g)$ and $(N^{n},h)$ are critical points of the energy functional 

$$E(\psi)=\frac{1}{2}\int_M|d\psi|^2dv_g,$$
\indent
where $\psi$ is a map from $M$ to $N$ and $dv_{g}$ denotes the volume element of $g$. The Euler-Lagrange equation of $E(\psi)$ is given by $\tau(\psi) = Trace\nabla d\psi = 0$, where $\tau(\psi)$ is the tension field of $\psi$, which vanishes precisely for harmonic maps.

In 1983, J. Eells and L. Lemaire \cite{EL} suggested to consider the problems associated to biharmonic maps which are a natural generalization of harmonic maps. A map $\psi$ is called {\it biharmonic} if it is a critical point of the bi-energy functional
$$E_2(\psi)=\frac{1}{2}\int_M|\tau(\psi)|^2dv_g,$$ 
on the space of smooth maps between two Riemannian manifolds. In \cite{Ji}, G.Y. Jiang studied the first and second variation formulas
of $E_{2}$ for which critical points are called biharmonic maps.  The Euler-Lagrange equation associated with this bi-energy functional is $\tau_2(\psi)=0$, where $\tau_2(\psi)$ is the so-called bi-tension field given by 
\begin{equation}\label{bitensionfield}
\tau_2(\psi)=\Delta\tau(\psi)-\tr\big(R^N(d\psi,\tau(\psi))d\psi\big).
\end{equation}
Here, $\Delta$ is the rough Laplacian acting on the sections of $\psi^{-1}(TN)$ given by $\Delta V=\tr(\nabla^2V)$ for any $V\in\Gamma(\psi^{-1}(TN))$  and $R^N$ is the curvature tensor of the target manifold $N$ defined as $R^N(X,Y)=[\nabla^N_X,\nabla^N_Y]-\nabla^N_{[X,Y]}$ for any $X,Y\in\Gamma(TN)$.

Over the past years, many geometers studied biharmonic
submanifolds and obtained a great variety of results in this domain (see \cite{BFO, CMO, CMO1, FLMO, FO, FOR, JI, LM, SH, MC, YO, YO1, Rot, RotAu, HU}, for instance). If the map $\psi: (M, g)\rightarrow (N, h)$ is an isometric immersion from a manifold $(M,g)$ into an ambient manifold $(N,h)$ then $M$ is called {\it biharmonic submanifold} of $N$. Since, it is obvious that any harmonic map is a biharmonic map, we will call {\it proper biharmonic submanifolds} the biharmonic submanifolds which are not harmonic, that is, minimal.

The main problem concerning biharmonic submanifold is the Chen's
Conjecture \cite{Ch2}:\linebreak

{\it ``\textbf{Biharmonic submanifolds of Euclidean spaces are the only submanifolds which are minimal ones.}''}\\

The Chen's biharmonic
conjecture is still an open problem, but lots of results on submanifolds of\linebreak Euclidean spaces provide affirmative partial solutions to the conjecture (see \cite{Chen,Ch3} and references\linebreak therein for an overview). On the other hand, the generalized Chen's conjecture replacing\linebreak Euclidean spaces by Riemannian manifolds of non-positive sectional curvature turns out to be false (see \cite{LO,OL} for counter-examples). Nevertheless, this generalized conjecture is true in \linebreak various situations and obtaining non-existence results in non-positive sectional curvature is still an interesting question. In \cite{RotAu}, authors gave two new contexts where such results hold.

In \cite{LU}, Lu gave a natural generalization of biharmonic maps and introduced $f$-biharmonic maps. He studied the first variation and calculated the $f$-biharmonic map equation as well as the equation for the $f$-biharmonic conformal maps between the same dimensional manfolds. Ou also studied $f$-biharmonic map and $f$-biharmonic submanifolds in \cite{YO2}, where he proved that an $f$-biharmonic map from a compact Riemannian manifold into a non-positively curved manifold with constant $f$-bienergy density is a harmonic map; any $f$-biharmonic function on a compact manifold is constant, and that the inversion about $S^{m}$ for $m\geq 3$ are proper $f$-biharmonic conformal diffeomorphisms. He also derived $f$-biharmonic submanifolds equation and proved that a surface in a manifold $(N^{n}, h)$ is an $f$-biharmonic surface if and only if it can be biharmonically conformally immersed into $(N^{n}, h)$. Further in \cite{YO3}, author characterize harmonic maps and minimal submanifolds by using the concept of $f$-biharmonic maps and obtained an improved equation for $f$-biharmonic hypersurfaces.
  
By definition, for a positive, well defined and $C^\infty$ differentiable function $f : M\rightarrow R$, \linebreak {\it $f$-biharmonic} maps are critical points of the $f$-bienergy functional for maps $\psi : (M,g)\rightarrow (N,h)$, between Riemannian manifolds, i.e.,
  
$$E_{2, f}(\psi)=\frac{1}{2}\int_M f|\tau(\psi)|^2dv_g.$$
Lu also obtained the corresponding Euler-Lagrange equation for $f$-biharmonic maps, i.e., 
\begin{equation}\label{relationfbiharmonic}
\tau_{2, f}(\psi)=f\tau_2(\psi) + (\Delta f)\tau(\psi) + 2\nabla^{\psi}_{\grad f}\tau(\psi) = 0.
\end{equation}
An $f$-biharmonic map is called a {\it proper $f$-biharmonic map} if it is neither a harmonic nor a biharmonic map. Also, we will call  {\it proper $f$-biharmonic submanifolds} a $f$-biharmonic submanifols which is neither minimal nor biharmonic.\\ \\
In \cite{OND}, the authors introduce another generalization of harmonic and biharmonic maps, namely, the $f$-harmonic and bi-$f$-harmonic maps. Given two Riemannian manifolds $(M^m,g)$ and $(N^n,h)$ and $f$ a smooth positive function over $M$,  they call bi-$f$-harmonic maps the critical points of the bi-$f$-energy functional for maps $\psi : (M,g)\rightarrow (N,h)$, between Riemannian manifolds:
$$E_f(\psi)=\frac12\int_Mf|\nabla \psi|^2dv_g.$$
The Euler-Lagrange equation is $\tau_f(\psi)=0$, where $\tau_f$ is the $f$-tension field defined by
$$\tau_f(\psi)=f\tau(\psi)+d\psi(\grad f).$$
Hence, the natural notion to consider is the bi-$f$-harmonicity given by the bi-$f$-energy functional
$$E_f^2(\psi)=\int_M|\tau_f(\psi)|^2dv_g.$$
Critical points of this functional are called bi-$f$-harmonic maps and are characterize by the following Euler-Lagrange equation
\begin{equation}\label{relationbifharmonic}
\tau^{2}_f(\psi)=fJ^{\psi}(\tau_f(\psi)) - \nabla^{\psi}_{\grad f}\tau_f(\psi) = 0,
\end{equation}
where  $J^{\psi}$ is the jacobi operator of the map defined by
 $$J^{\psi}(X) = -[Tr_{g}\nabla^{\psi}\nabla^{\psi}X - \nabla^{\psi}_{\nabla^{M}}X - R^{N}(d\psi, X)d\psi].$$
Obviously, $f$-harmonic maps are bi-f-harmonic maps, hence, we will call {\it proper $f$-biharmonic map} a $f$-biharmonic map which is not $f$-harmonic. However, we want to point out that there is no link between the notion of $f$-biharmonic and bi-$f$-harmonic maps. In particular, there is no reason for a $f$-harmonic maps to bi $f$-biharmonic.\\ \\
\indent 
In the present paper, we will focus here on  $f$-biharmonic submanifolds and bi-$f$-harmonic submanifolds of both (generalized) complex space forms and generalized Sasakian space forms. The paper is organized as follows: \\
In Section \ref{sec2}, we recall the basics of genaralized complex and Sasakian space forms as well as their submanifolds. Section \ref{sec3} is devoted to the study of $f$-biharmonic submanifolds. For both classes of ambient spaces, we first give the general necessary and sufficient condition for submanifolds to be $f$-biharmonic. Then, we focus of many particular cases and obtain some non-existence results. Finally, Section \ref{sec4} is devoted to bi-$f$-harmonic submanifolds. At first, since the notion of bi-$f$-harmonic submanifold almost has not been studied, we give a general characterization  of bi-$f$-harmonic submanifold in any ambient space. Then, we apply this general result to the case of generalized complex and Sasakian space forms. \\ \\
\indent

\section{\textbf{Preliminaries}}\label{sec2}
\subsection{Generalized complex space forms and their submanifolds} A Hermitian manifold $(N,g,J)$ with constant sectional holomorphic curvature $4c$ is called a complex space form. We\linebreak denote by $M^n_{\CC}(4c)$ be the simply connected complex $n$-dimensional complex space form of \linebreak constant holomorphic sectional curvature $4c$. The curvature tensor $R$ of $M^n_{\CC}(4c)$ is given by
\begin{center}
$R^\CC(X,Y)Z
= c\{g(Y,Z)X - g(X,Z)Y + g(Z,JY)JX - g(Z,JX)JY + 2g(X,JY)JZ\}$,
\end{center} 
for $X, Y, Z\in \Gamma(TM^n_{\CC}(4c))$, where $<\cdot, \cdot>$ is the Riemannian metric on $M^n_{\CC}(4c)$ and $J$ is the
almost complex structure of $M^n_{\CC}(4c)$. The complex space from $M^n_{\CC}(4c)$ is the complex projective space $\CC P^{n}(4c)$, the complex Euclidean space $\CC^{n}$ or the complex hyperbolic
space $\CC H^{n}(4c)$ according to $c > 0$, $c = 0$ or $c < 0$.\\ \\
Now, we consider a natural generalization of complex space forms, namely the generalized\linebreak complex space forms. After defining them, we will give some basic information about generalized complex space forms and their submanifolds. Generalized complex space forms form a particular class of Hermitian manifolds which has not been intensively studied. In 1981, Tricelli and Vanhecke \cite{TV} introduced the following generalization of the complex space forms ($\CC^n$, $\CC P^n$ and $\CC H^n$). Let $(N^{2n},g,J)$ be an almost Hermitian manifold. We denote the generalized curvature tensors by $R_1$ and $R_2$ which is defined as
$$R_1(X,Y)Z=g(Y,Z)X-g(X,Z)Y,$$
$$R_2(X,Y)Z=g(JY,Z)JX-g(JX,Z)JY+2g(JY,X)JZ, \hspace{.2cm}\forall\hspace{.2cm} X,Y,Z\in \Gamma(TN).$$
The manifold $(N,g,J)$ is called {\it generalized complex space form} if its curvature tensor $R$ has the following form
$$R=\alpha R_1+\beta R_2,$$
where $\alpha$ and $\beta$ are smooth functions on $N$. The terminology comes obviously from the fact that complex space forms satisfy this property with constants $\alpha=\beta$ . \\
\indent
In the same paper \cite{TV}, Tricelli and Vanhecke showed that if $N$ is of (real) dimension $2n\geq6$, then $(N,g,J)$ is a complex space form. They also showed that $\alpha+\beta$ is necessarily constant. This implies that $\alpha=\beta$ are constants in dimension $2n\geq6$, but this is not the case in dimension $4$. Hence, the notion of generalized complex space form is of interest only in dimension $4$.  Further, Olszak \cite {Ols} constructed examples in dimension $4$ with $\alpha$ and $\beta$ non-constant. These examples are obtained by conformal deformation of B\"ochner flat K\"ahlerian manifolds of non constant scalar curvature.  Examples of B\"ochner flat K\"ahlerian manifolds can be found in \cite{Der}. From now on, we will denote by $N(\alpha,\beta)$ a (4-dimensional) generalized complex space form with curvature given by $R=\alpha R_1+\beta R_2$. Note that these spaces are Einstein, with constant scalar curvature equal to $12(\alpha+\beta)$. Of course, they are not K\"ahlerian because if they were, they would be complex space forms. \\ \\
\indent
Now, let $M$ be a submanifold of the (generalized) complex space form $M^n_{\CC}(4c)$ or $N(\alpha,\beta)$. The almost complex structure $J$ on $M^n_{\CC}(4c)$ (or $N(\alpha,\beta)$) induces the existence of four operators on $M$, namely
$$j:TM\lgra TM,\ k:TM\lgra NM,\l:NM\lgra TM\ \text{and}\ m:NM\lgra NM ,$$
defined for all $X\in TM$ and all $\xi\in NM$ by
\begin{eqnarray}\label{relationfhst}
JX=jX+kX\quad\text{and}\quad
J\xi=l\xi+m\xi.
\end{eqnarray} 
Since $J$ is an almost complex structure, it satisfies $J^2=-Id$ and for $X,Y$ tangent to $M^n_{\CC}(4c)$ (or $N(\alpha,\beta)$), we have $g(JX,Y)=-g(X,JY)$. Then, we deduce that the operators $j,k,l,m$ satisfy the following relations
\begin{align}
&j^2X+lkX=-X,& \label{relation1.1}\\
&m^2\xi+kl\xi=-\xi,& \label{relation1.2}\\
\label{relation1.3}
&jl\xi+lm\xi=0,&\\
\label{relation1.4}
&kjX+mkX=0,&\\
\label{relation1.5}
&g(kX,\xi)=-g(X,l\xi),&
\end{align}
for all $X\in\Gamma(TM)$ and all $\xi\in\Gamma(NM)$. Moreover $j$ and $m$ are skew-symmetric.

\subsection{Generalized Sasakian space forms and their submanifolds}
Now, we give some \linebreak recalls about almost contact metric manifolds and generalized Sasakian space forms. For more details, one can refer to (\cite{ABC,Bla,YK}) for instance. A Riemannian manifold $\widetilde{M}$ of odd dimension is said almost contact if there exists globally over $\widetilde{M}$, a  vector field $\xi$, a $1$-form $\eta$ and a field of $(1,1)$-tensor $\phi$ satisfying the following conditions:
\begin{equation} \eta(\xi)=1\quad\text{and}\quad \phi^2=-Id+\eta\otimes\xi.
\end{equation}
Remark that this implies $\phi\xi=0$ and $\eta\circ\phi=0$. The manifold $\widetilde{M}$ can be endowed with a Riemannian metric $\widetilde{g}$ satisfying
\begin{equation}
\widetilde{g}(\phi X,\phi Y)=\widetilde{g}(X,Y)-\eta(X)\eta(Y)\quad\text{and}\quad \eta(X)=\widetilde{g}(X,\xi),
\end{equation}
for any vector fields $X,Y$ tangent to $\widetilde{M}$. Then, we say that $(\widetilde{M},\widetilde{g},\xi,\eta,\phi)$ is an almost contact metric manifold. Three class of this family are of particular interest, namely, the Sasakian, Kenmotsu and cosymplectic manifolds. We will give some recalls about them. \\ \\
\indent
First, we introduce the fundamental $2$-form (also called Sasaki $2$-form) $\Omega$ defined for $X,Y\in\Gamma(TM)$ by 
$$\Omega(X,Y)=\widetilde{g}(X,\phi Y).$$
We consider also $N_{\phi}$, the Nijenhuis tensor defined by
$$N_{\phi}(X,Y)=[\phi X,\phi Y]-\phi[\phi X,Y]-\phi[X,\phi Y]+\phi^2[X,Y],$$
for any vector fields $X,Y$. An almost contact metric manifold is said normal if and only if the Nijenhuis tensor $N_{\phi}$ satisfies
$$N_{\phi}+2d\eta\otimes\xi=0.$$ 
An almost contact metric manifold is said {\it Sasakian manifold} if and only if it is normal and  $d\eta=\Omega$. This is equivalent to
\begin{equation}
(\nabla_X\phi)Y=\widetilde{g}(X,Y)\xi-\eta(Y)X, \hspace{.2 cm}\forall \hspace{.2 cm}X, Y\in\Gamma(\widetilde{M}).
\end{equation}
It also implies that 
\begin{equation}
\nabla_X\xi=-\phi(X).
\end{equation}
An almost contact metric manifold is said {\it Kenmotsu manifold} if and only if $d\eta=0$ and $d\Omega=2\eta\wedge\Omega$. Equivalently, this means
\begin{equation}
(\nabla_X\phi)Y=-\eta(Y)\phi X-g(X,\phi Y)\xi,
\end{equation}
for any $X$ and $Y$. Hence, we also have 
\begin{equation}
\nabla_X\xi=X-\eta(X)\xi.
\end{equation}
Finally, an almost contact metric manifold is said {\it cosymplectic manifold} if and only if $d\eta=0$ and $d\Omega=0$, or equivalently 
\begin{equation}
\nabla\phi=0,
\end{equation}
and in this case, we have
\begin{equation}
\nabla\xi=0.
\end{equation}
The $\phi$-sectional curvature of  an almost contact metric manifold is  defined as the sectional \linebreak curvature on the $2$-planes $\{X,\phi X\}$. When the $\phi$-sectional curvature is constant, we say that the manifold is a space form (Sasakian, Kenmotsu or cosymplectic in each of the three cases above). It is well known that the $\phi$-sectional curvature determines entirely the curvature of the manifold. When the $\phi$-sectional curvature is constant, the curvature tensor is expressed explicitely. Let $R_1^{\star}$, $R_2^{\star}$ and $R_3^{\star}$ be the generalized curvature tensors defined by
\begin{equation}
R_1^{\star}(X,Y)Z=\widetilde{g}(Y,Z)X-\widetilde{g}(X,Z)Y,
\end{equation}
\begin{equation}
R_2^{\star}(X,Y)Z=\eta(X)\eta(Z)Y-\eta(Y)\eta(Z)X+\widetilde{g}(X,Z)\eta(Y)\xi-\widetilde{g}(Y,Z)\eta(X)\xi 
\end{equation} and 
\begin{equation}
R_3^{\star}(X,Y)Z=\Omega(Z,Y)\phi X-\Omega(Z,X)\phi Y+2\Omega(X,Y)\phi Z. 
\end{equation}
For the three cases we are interested in, the curvature of a space form of constant $\phi$-sectional curvature $c$ is given by

\begin{itemize}
\item Sasaki: $R^{\star}=\frac{c+3}{4}R_1^{\star}+\frac{c-1}{4}R_2^{\star}+\frac{c-1}{4}R_3^{\star}.$\\
\item Kenmotsu: $R^{\star}=\frac{c-3}{4}R_1^{\star}+\frac{c+1}{4}R_2^{\star}+\frac{c+1}{4}R_3^{\star}.$\\
\item Cosymplectic: $R^{\star}=\frac{c}{4}R_1^{\star}+\frac{c}{4}R_2^{\star}+\frac{c}{4}R_3^{\star}.$
\end{itemize}
In the sequel, for more clarity, we will denote the Sasakian ({\it resp.}  Kenmotsu, cosymplectic) space form of constant $\phi$-sectional curvature $c$ by $\widetilde{M}_S(c)$ ({\it resp. } $\widetilde{M}_K(c)$, $\widetilde{M}_C(c)$). These space forms appear as particular cases of the so-called generalized Sasakian space forms, introduced by Alegre, Blair and Carriazo in \cite{ABC}. A generalized Sasakian space form, denoted by $\widetilde{M}(f_1,f_2,f_3)$, is a contact metric manifold with curvature tensor of the form
\begin{equation}\label{CurvatureGSasakian}
f_1R_1^{\star}+f_2R_2^{\star}+f_3R_3^{\star},
\end{equation}
where $f_1$, $f_2$ and $f_3$ are real functions on the manifold. The most simple examples of generalized Sasakian space forms are the  warped products of the real line by a complex space form or a generalized complex space forms.  Their conformal deformations as well as their so-called $\mathcal{D}$-homothetic deformations are also generalized Sasakian space forms (see \cite{ABC}). Other examples can be found in \cite{AC}.\\ \\
Now, let $(M,g)$ be a submanifold of an almost contact metric manifold $(\widetilde{M},\widetilde{g},\xi,\eta,\phi)$. The field of tensors $\phi$ induces on $M$, the existence of the following four operators:
$$P:TM\lgra TM,\ N:TM\lgra NM,\ t:NM\lgra TM\ \text{and}\ s:NM\lgra NM ,$$
defined for any $X\in TM$ and $\nu\in NM$. Now, we have
\begin{eqnarray}\label{relationfhst}
\phi X=PX+NX\quad\text{and}\quad
\phi\nu=s\nu+t\nu,
\end{eqnarray} 
where $PX$ and $NX$ are tangential and normal components of $\phi X$, respectively, whereas $t\nu$ and $s\nu$ are the tangential and normal components of $\phi\nu$, respectively. A submanifold $M$ is said invariant ({\it resp.} anti-invariant) if $N$ ({\it resp.} $P$) vanishes identically. In \cite{Lot}, Lotta shows that if the vector field $\xi$ is normal to $M$, then $M$ is anti-invariant.  
\section{$f$-Biharmonic submanifolds}\label{sec3}

\subsection{\textbf{$f$-Biharmonic submanifolds of generalized complex space forms}}\label{sec31}
At first, we will calculate necessary and sufficient condition of $f$-biharmonic submanifold of generalized complex space forms and then we make a exposition about the results which could characterize these type of submanifolds.
\begin{thm}\label{thm1}
Let $M^{p}$, $p<4$ be a submanifold of the generalized complex space form $N(\alpha,\beta)$ with second fundamental form $B$, shape operator $A$, mean curvature $H$ and a positive $C^{\infty}$-differentiable function $f$ on $M$. Then $M$ is $f$-biharmonic submanifold of $N(\alpha,\beta)$ if and only if the  following two equations are satisfied
\begin{enumerate}
\item $$
-\Delta^{\perp}H+\tr \left(B(\cdot,A_H\cdot)\right)-p\alpha H+3\beta klH + \frac{\Delta f}{f}H + 2\nabla^{\perp}_{\grad(\ln f)}H=0,$$
\item $$ \frac{p}{2}{\rm \grad}|H|^2-2A_{H} \grad(\ln f)+2\tr \left(A_{\nabla^{\perp}H}(\cdot)\right)+6\beta jlH=0.$$
\end{enumerate}
\end{thm}
\noindent
{\bf Proof:} 
It is a classic fact that the tension field of the isometric immersion $\psi$ is given by
\begin{equation}\label{tensionfield}
\tau(\psi) = \tr\nabla d\psi = \tr B = p H.
\end{equation}

Using equation (\ref{tensionfield}) in equation (\ref{bitensionfield}), we have
\begin{equation}\label{eqtau2}
\tau_2(\psi)=p\Delta H-\tr\big(R^N(d\psi,pH)d\psi\big).\end{equation}
Moreover, we recall that, by some classical and straightforward computations, we have
$$\Delta H=\frac{p}{2}\grad|H|^2+\tr \left(B(\cdot,A_H\cdot)\right)+2\tr \left(A_{\nabla^{\perp}H}(\cdot)\right)+\Delta^{\perp}H.$$
Reporting this into \eqref{eqtau2}, we get
\begin{equation}\label{eqbiharmonic}
\tau_2(\psi)=-\Delta^{\perp}H+\tr \left(B(\cdot,A_H\cdot)\right)+
\frac{p}{2}{\rm grad}|H|^2+2\tr \left(A_{\nabla^{\perp}H}(\cdot)\right)+2\tr\left(R^N(\cdot,H)\cdot\right).
\end{equation}

Now, the curvature tensor of generalized complex space form, $N(\alpha,\beta)$, is given by
\begin{eqnarray*}
\tr\left(R(\cdot,H)\cdot\right)&=&\alpha\tr\left(R_1(\cdot,H)\cdot\right)+\beta\tr\left(R_2(\cdot,H)\cdot\right).
\end{eqnarray*}
Let $\{e_1,\cdots,e_p\}$ be a local orthonormal frame of $TM$. Then, we have
\begin{center}
$\tr\left(R(\cdot,H)\cdot\right)= \alpha\displaystyle\sum_{i=1}^pR_1(e_i,H)e_i + \beta\sum_{i=1}^pR_2(e_i,H)e_i $
\end{center}
or,
\begin{eqnarray*}
\tr\left(R(\cdot,H)\cdot\right)=\alpha\sum_{i=1}^p\left[g(H,e_i)e_i-g(e_i,e_i)H\right] \\+ \beta\sum_{i=1}^p\left[g(JH,e_i)Je_i-g(Je_i,e_i)JH+2g(JH,e_i)Je_i\right].
\end{eqnarray*}
or, 
\begin{eqnarray}\label{curvfbiharm}
\tr\left(R(\cdot,H)\cdot\right)=\alpha(-pH) + \beta(3jlH + 3klH). 
\end{eqnarray}

From equation (\ref{relationfbiharmonic}), $M$ is $f$-biharmonic if and only if  
\begin{center}
$f\tau_2(\psi) + \Delta f\tau(\psi) + 2\nabla^{\psi}_{\grad f}\tau(\psi) = 0,$
\end{center}
which is equivalent to
\begin{equation}\label{fbiharmoniccondition}
\tau_2(\psi) + p\frac{\Delta f}{f}H + 2p(-A_{H} \grad(\ln f) + \nabla^{\perp}_{\grad(\ln f)}H) = 0.
\end{equation}
Now, using equations \eqref{eqbiharmonic} and \eqref{curvfbiharm} in equation \eqref{fbiharmoniccondition} and considering that $jlH$ is tangent and $klH$ is normal, we get the statement of the theorem by identification of tangent and normal parts.  \hfill $\square$
\begin{cor}\label{cor1}
Let $M^{p}$,  $p\leqslant 2n$, be a submanifold of the complex space form $M^n_{\CC}(4c)$ of complex dimension $n$ and constant holomorphic sectional curvature $4c$, with second fundamental form $B$, shape operator $A$, mean curvature $H$ and a positive $C^{\infty}$-differentiable function $f$ on $M$. Then $M$ is $f$-biharmonic submanifold of $M^n_{\CC}(4c)$ if and only if the following two equations are satisfied

\begin{enumerate}
\item $$
-\Delta^{\perp}H+\tr \left(B(\cdot,A_H\cdot)\right)-pc H+3c klH + \frac{\Delta f}{f}H + 2\nabla^{\perp}_{\grad(\ln f)}H=0,$$

\item $$\frac{p}{2}{\rm grad}|H|^2-2A_{H} \grad(\ln f)+2\tr \left(A_{\nabla^{\perp}H}(\cdot)\right)+6c jlH=0.
$$
\end{enumerate}
\end{cor}
{\bf Proof:} For complex space forms the computations are essentially the same as for the generalized complex space forms with the only differences that $\alpha=\beta=c$ and dimension is not necessarily equal to 4. 
\hfill $\square$\\ \\
In the sequel, we will state many results for biharmonic subamnifolds of the generalized complex space forms $N(\alpha,\beta)$. They have of course analogue for the complex space forms but for a sake of briefness, we do not write then since the results are the same with $\alpha=\beta=c$.  
Assuming particular cases such as hypersurfaces, Lagrangian or complex surfaces and curves of generalized complex space form $N(\alpha,\beta)$, we have the following conclusion.
\begin{cor}\label{cor1}
Let $M^{p}$, $p<4$ be a submanifold of the generalized complex space form $N(\alpha,\beta)$ with second fundamental form $B$, shape operator $A$, mean curvature $H$ and a positive $C^{\infty}$-differentiable function $f$ on $M$.
\begin{enumerate}
\item If $M$ is a hypersurface then $M$ is $f$-biharmonic if and only if

 $$-\Delta^{\perp}H + \frac{\Delta f}{f}H + 2\nabla^{\perp}_{\grad(\ln f)}H + \tr \left(B(\cdot,A_H\cdot)\right)-3(\alpha+\beta) H=0,$$
and
$$\frac{3}{2}{\rm \grad}|H|^2 - 2A_{H}grad(\ln f)+2\tr \left(A_{\nabla^{\perp}H}(\cdot)\right)=0.$$

\item If $M$ is a complex surface then $M$ is $f$-biharmonic if and only if
$$-\Delta^{\perp}H + \frac{\Delta f}{f}H + 2\nabla^{\perp}_{\grad(\ln f)}H + \tr \left(B(\cdot,A_H\cdot)\right)-2\alpha H=0,$$
and
$${\rm \grad}|H|^2-2A_{H} \grad(\ln f)+2\tr \left(A_{\nabla^{\perp}H}(\cdot)\right)=0.$$
\item If $M$ is a Lagrangian surface then $M$ is $f$-biharmonic if and only if
$$-\Delta^{\perp}H + \frac{\Delta f}{f}H + 2\nabla^{\perp}_{\grad(\ln f)}H + \tr \left(B(\cdot,A_H\cdot)\right)-2\alpha H-3\beta H=0,$$
and
$${\rm \grad}|H|^2 - 2A_{H}\grad(\ln f) +2\tr \left(A_{\nabla^{\perp}H}(\cdot)\right)=0.
$$
\item If $M$ is a curve then $M$ is $f$-biharmonic if and only if
$$-\Delta^{\perp}H +\frac{\Delta f}{f}H + 2\nabla^{\perp}_{\grad(\ln f)}H +\tr\left( B(\cdot,A_H\cdot)\right)-\alpha H-3\beta (H+m^2H)=0,$$
and
$$ \frac{1}{2}{\rm \grad}|H|^2 - 2A_{H}\grad(\ln f) +2\tr \left(A_{\nabla^{\perp}H}(\cdot)\right)=0.$$
\end{enumerate}
\end{cor}
{\bf Proof:} The proof is a consequence of Theorem \ref{thm1} using the facts that
\begin{enumerate}
\item if $M$ is a hypersurface, then $m=0$ and so $jlH=0$, $kjH=0$ and $klH=-H$,
\item if $M$ is a complex surface then $k=0$ and $l=0$,
\item if $M$ is a Lagrangian surface, then $j=0$, $m=0$,
\item if $M$ is a curve, then $j=0$.
\end{enumerate}
\hfill$\square$

\begin{remark}
It is a well known fact that any complex submanifold of a K\"ahler manifold is \linebreak necessarily minimal. But as mentioned above, the generalized space forms $N(\alpha,\beta)$ are not K\"ahlerian unless there are the complex projective plane or the complex hyperbolic plane. Hence, considering $f$-biharmonic surfaces into $N(\alpha,\beta)$ is of real interest, since they are not necessarily minimal.
\end{remark}

Similarly, if we assume mean curvature vector $H$ as parallel vector then for curves and complex or Lagranian surfaces, we obtain the following corollaries.   
\begin{cor}\label{corlag}
Let $M^{p}$, $p<4$ be a submanifold of the generalized complex space form $N(\alpha,\beta)$ with second fundamental form $B$, shape operator $A$, mean curvature $H$ and a positive $C^{\infty}$-differentiable function $f$ on $M$.
\begin{enumerate}
\item If $M$ be a Lagrangian surface of $N(\alpha,\beta)$ with parallel mean curvature then $M$ is $f$-biharmonic if and only if
\begin{center}
$\tr \left(B(\cdot,A_H\cdot)\right)= 2\alpha H+3\beta H - \frac{\Delta f}{f}H,\hspace{.2cm} and\hspace{.2cm}
A_{H}\grad f=0$.
\end{center}

\item If $M$ be a complex surface of $N(\alpha,\beta)$ with parallel mean curvature then $M$ is $f$-biharmonic if and only if
\begin{center}
$\tr\left(B(\cdot,A_H\cdot)\right)=2\alpha H -\frac{\Delta f}{f}H\hspace{.2cm}  \textit{and} \hspace{.2cm}
A_{H}\grad f=0$.
\end{center}
\item If $M$ is a curve in $N(\alpha,\beta)$ with parallel mean curvature then $M$ is $f$-biharmonic if and only if
\begin{center}
$\tr\left( B(\cdot,A_H\cdot)\right)= \alpha H+3\beta (H+m^2H)-\frac{\Delta f}{f}H,\hspace{.2cm} and\hspace{.2cm}
A_{H}\grad f=0.$
\end{center}
\end{enumerate}
\end{cor}
{\it Proof:} Since $M$ has parallel mean curvature so that the terms $\Delta^{\perp}H$, $\nabla^{\perp}_{\grad f}H$, $\grad|H|^2$ and $\tr(A_{\nabla^{\perp}_{\cdot}H\cdot})$ vanish and we obtain immediately the result from the previous Corollary.
\hfill$\square$\\ 
\begin{remark}
Note that for the last two results there is no analogue for complex subamnifolds of $M^n_{\CC}(4c)$ since they are in fact minimal.
\end{remark}
Further, for  constant mean curvature hypersurfaces in $N(\alpha,\beta)$, we have the following result.

\begin{prop}\label{propB}
$(1)$ Let $M^{3}$ be a hypersurface of the generalized complex space form $N(\alpha,\beta)$ with second fundamental form $B$, non zero constant mean curvature $H$ and $f$ a positive $C^{\infty}$-differentiable function on $M$. Then $M$ is $f$ biharmonic if and only if
$$|B|^2=3(\alpha+\beta)- \frac{\Delta f}{f}\quad\text{and}\quad A\,\grad f=0$$
or equivalently, $M$ is proper $f$-biharmonic if and only if the scalar curvature of $M$ satisfies
$$\scal_M=3(\alpha+\beta)+9H^2+\frac{\Delta f}{f}\quad\text{and}\quad A\,\grad f=0.$$
$(2)$ There exists no proper $f$-biharmonic hypersurfaces with constant mean curvature and \linebreak constant scalar curvature.
\end{prop}

\noindent
{\bf Proof:} For the first point, since $M$ is a hypersurface, by Corollary \ref{cor1}, $M$ is $f$-biharmonic if and only if
$$\left\{
\begin{array}{l}
-\Delta^{\perp}H + \frac{\Delta f}{f}H + 2\nabla^{\perp}_{\grad(\ln f)}H + \tr \left(B(\cdot,A_H\cdot)\right)-3(\alpha+\beta) H=0,\\ \\
\frac{3}{2}{\rm \grad}|H|^2 - 2A_{H}\grad(\ln f)+2\tr \left(A_{\nabla^{\perp}H}(\cdot)\right)=0.
\end{array}
\right.$$
Since $M$ has constant mean curvature, the above equation reduces to
$$\left\{
\begin{array}{l}
\tr \left(B(\cdot,A_H\cdot)\right)= 3(\alpha+\beta) H - \frac{\Delta f}{f}H,\\ \\
A_{H}\grad(\ln f)=0.
\end{array}
\right.$$

Using condition $A_H=HA$ for hypersurfaces, we get
\begin{center}
$\tr \Big(B(\cdot,A_H(\cdot))\Big)=H\tr \Big(B(\cdot,A(\cdot))\Big)=H|B|^2.$
\end{center}
Reporting this result in first equation of the above condition and from the assumption that $H$ is a non-zero constant, we get the desired identity $|B|^2=3(\alpha+\beta) - \frac{\Delta f}{f}$.\\
For the second equivalence, by the Gauss equation, we have
$$
\scal_M=\sum_{i,j=1}^3g\left( R^N(e_i,e_j)e_j,e_i\right)-|B|^2+9H^2,
$$
where $\{e_1,e_2,e_3\}$ is a local orthonormal frame of $M$. From the expression of the curvature tensor of $N(\alpha,\beta)$, we get
$$
\scal_M=6(\alpha+\beta)-||B||^2+9H^2.$$
Moreover, since $\grad (\ln f)=\frac{1}{f}\grad f$ and $A_H=HA$ with $H$ is a non-zero constant, then $A_{H}\grad(\ln f)=0$ reduces to $A\,\grad f=0$. 

Hence, we deduce that $M$ is proper $f$-biharmonic if and only if $|B|^2=3(\alpha+\beta)-(\frac{\Delta f}{f})$ and $A\,\grad f=0$, that is, if and only if $\scal_M=3(\alpha+\beta)+9H^2 +\frac{\Delta f}{f}$ and $A\,\grad f=0$.\\ \\
Now, for the second point, if $M$ is a hypersurface with constant mean curvature and constant scalar curvature, then by the first point, if $M$ is $f$-biharmonic then
$$\scal_M=3(\alpha+\beta)+9H^2+\frac{\Delta f}{f}.$$
As we have already mentioned, $\alpha+\beta$ is constant, hence, since $H$ and $\scal_{M}$ are constant, then $\frac{\Delta f}{f}$ is constant, that is, $f$ is an eigenvalue of the Laplacian. But $f$ is a positive function, so the only possibility is that $f$ is a positive constant and $M$ is biharmonic. This concludes the proof of the second point.
\hfill$\square$\\ \\
Now, we give this proposition which give an estimate of the mean curvature for a $f$-biharmonic Lagrangian surface. 
\begin{prop}\label{proplag}
Let $M^{2}$ be a Lagrangian surface of the generalized complex space form $N(\alpha,\beta)$ with second fundamental form $B$, shape operator $A$, non-zero constant mean curvature $H$ and a positive $C^{\infty}$-differentiable function $f$ on $M$. 
\begin{enumerate} 
\item If $\inf_M\left(2\alpha+3\beta-\frac{\Delta f}{f}\right)$ is non-positive then $M$ is not $f$-biharmonic. \\
\item
If $\inf_M\left(2\alpha+3\beta-\frac{\Delta f}{f}\right)$ is positive and $M$ is proper $f$-biharmonic then $$0<|H|^2 \leqslant \inf_{M}\left(\frac{2\alpha+3\beta-\frac{\Delta f}{f}}{2}\right).$$ 
\end{enumerate}
\end{prop}
{\bf Proof:} Assume that $M$ is a $f$-biharmonic Lagrangian surface of $N(\alpha,\beta)$, considering third assertion of Corollary \ref{cor1}, we have

$$\left\{
\begin{array}{l}
-\Delta^{\perp}H + \frac{\Delta f}{f}H + \frac{2}{f}\nabla^{\perp}_{\grad f}H + \tr \left(B(\cdot,A_H\cdot)\right)-2\alpha H-3\beta H=0,\\ \\
{\rm \grad}|H|^2 - \frac{2}{f}A_{H}\grad f +2\tr \left(A_{\nabla^{\perp}H}(\cdot)\right)=0.
\end{array}
\right.$$
Hence, by taking the scalar product with $H$ and taking the assumption that mean curvatutre $H\neq 0$, i.e., $|H|$ is constant,  from the first part of the above equation, we have
\begin{center}
$-<\Delta^{\perp}H, H> + \frac{2}{f}<\nabla^{\perp}_{\grad f}H,H>+  |A_H|^2-\left(\frac{\Delta f}{f}-2\alpha -3\beta\right) <H, H>=0.$
\end{center}
This equation implies that
$$-\left\langle\Delta^{\perp}H,H\right\rangle=\left(2\alpha+3\beta-\frac{\Delta f}{f}\right)|H|^2-|A_H|^2,$$
where we have used that $<\nabla^{\perp}_{\grad f}H,H>=0$ since $|H|$ is constant. Now, with the help of the Bochner formula, we get
$$\left(2\alpha+3\beta-\frac{\Delta f}{f}\right)|H|^2=|A_H|^2+|\nabla^{\perp}H|^2 .$$
Now, using Cauchy-Schwarz inequality, i.e., $|A_H|^2\geqslant 2|H|^4$ in the above equation, we have 
\begin{equation}\label{inequality}
\left(2\alpha+3\beta-\frac{\Delta f}{f}\right)|H|^2  \geqslant 2|H|^4+|\nabla^{\perp}H|^2 \geqslant 2|H|^4 .
\end{equation}
So, we have 
$0<|H|^2\leqslant \inf_M\left(\frac{2\alpha+3\beta-\frac{\Delta f}{f}}{2}\right)$ because $|H|$ is a non-zero constant. This is only possible if the function $2\alpha+3\beta-\frac{\Delta f}{f}$ has a positive infimum. This concludes the proof.
\hfill $\square$\\ \\
Now, we have similar result for complex surfaces.
\begin{prop}\label{propcomp}
Let $\psi: M^{2}\rightarrow N(\alpha,\beta)$ be a complex surface of generalized complex space form $N(\alpha,\beta)$ with second fundamental form $B$, shape operator $A$, mean curvature $H$ and a positive $C^{\infty}$-differentiable function $f$ on $M$. 
\begin{enumerate} 
\item If $\inf_M\left(2\alpha-\frac{\Delta f}{f}\right)$ is non-positive then $M$ is not $f$-biharmonic. \\
\item
If $\inf_M\left(2\alpha-\frac{\Delta f}{f}\right)$ is positive and $M$ is proper $f$-biharmonic then $$0<|H|^2 \leqslant \inf_{M}\left(\frac{2\alpha-\frac{\Delta f}{f}}{2}\right).$$  

\end{enumerate}
\end{prop}
{\bf Proof:} Let $M$ be a $f$-biharmonic complex surface of $N(\alpha,\beta)$ with non-zero constant mean curvature. Then, by the second assertion of Corollary \ref{cor1}, we have
\begin{center}
$-\Delta^{\perp}H + \frac{\Delta f}{f}H + \tr \left(B(\cdot,A_H\cdot)\right)-2\alpha H=0,\hspace{.2cm}and\hspace{.2cm}
A_{H}\grad f=0.$
\end{center}
Replacing $2\alpha+3\beta$ by $2\alpha$ in the proof of Proposition \ref{proplag}, we have the required result.
\hfill $\square$\\ \\

\subsection{$f$-Biharmonic submanifolds of generalized Sasakian space forms}\label{sec32}
Now, we consider $f$-biharmonic submanifolds of generalized Sasakian space forms and give the following theorm for its characterization.
\begin{thm}\label{thm2}
Let $M^{p}$ be a submanifold of a generalized Sasakian space form $\widetilde{M}(f_1,f_2,f_3)$, with second fundamental form $B$, shape operator $A$, mean curvature $H$ and a positive $C^{\infty}$-differentiable function $f$ on $M$. Then $M$ is $f$-biharmonic submanifold of $\widetilde{M}(f_1,f_2,f_3)$ if and only if the following two equations are satisfied
\begin{eqnarray*}
-\Delta^{\perp}H+\tr B(\cdot,A_H\cdot) + \frac{\Delta f}{f}H +  2\nabla^{\perp}_{\grad(\ln f)}H=pf_1H-f_2|\xi^{\top}|^2H-pf_2\eta(H)\xi^{\perp}-3f_3NsH
\end{eqnarray*}
and
$$\frac{p}{2}\grad|H|^2+2\tr A_{\nabla^{\perp}H}(\cdot)-2A_{H} \grad(\ln f)=-2f_2(p-1)\eta(H)\xi^{\top}-6f_3PsH.$$
\end{thm}
{\bf Proof:} At first, we calculate the curvature tensor of generalized Sasakian space form $\widetilde{M}(f_1,f_2,f_3)$. From equation \eqref{CurvatureGSasakian}, we have

\begin{eqnarray*}
R^{\star}(X, Y)Z &=& f_{1}R_{1}^{\star}(X, Y)Z + f_{2}R_{1}^{\star}(X, Y)Z + f_{3}R_{2}^{\star}(X, Y)Z \\
&=& f_{1}\{\tilde{g}(Y, Z)X - \tilde{g}(X, Z)Y\}\\
&+& f_{2}\{\eta(X)\eta(Z)Y - \eta(Y)\eta(Z)X + \tilde{g}(X, Z)\eta(Y)\xi - \tilde{g}(Y, Z)\eta(X)\xi \} \\ &+& f_{3}\{\tilde{g}(X, \phi Z)\phi Y - \tilde{g}(Y, \phi Z)\phi X + 2\tilde{g}(X, \phi Y)\phi Z\}.
\end{eqnarray*}

Let us consider $\{e_{1}, e_{2}, ..., e_{p}\}$ an orthogonal basis of the tangent space of $M$. Then, we have

\begin{eqnarray*}
R^{\star}(e_{i}, H)e_{i} &=& f_{1}\{\tilde{g}(H, e_{i})e_{i} - \tilde{g}(e_{i}, e_{i})H\}
+ f_{2}\{\eta(e_{i})\eta(e_{i})H - \eta(H)\eta(e_{i})e_{i} + \tilde{g}(e_{i}, e_{i})\eta(H)\xi \} \\ &+& f_{3}\{\tilde{g}(e_{i}, \phi e_{i})\phi H - \tilde{g}(H, \phi e_{i})\phi e_{i} + 2\tilde{g}(e_{i}, \phi H)\phi e_{i}\}.
\end{eqnarray*}

Taking the trace and using \eqref{relationfhst} in the above equation, we get 

\begin{eqnarray*}
tr \big(R^{\star}(\cdot, H)\cdot\big) &=&  - f_{1}pH 
+ f_{2}\sum_{i}\{\eta(e_{i})^{2}H - \eta(H)\eta(e_{i})e_{i} + |e_{i}|^{2}\eta(H)\xi \} \\ &+& f_{3}\sum_{i}\{tr(P)\phi H - \tilde{g}(H, N e_{i})\phi e_{i} + 2\tilde{g}(e_{i}, s H)\phi e_{i}\}
\\
&=&- f_{1}pH 
+ f_{2}\{|\xi^{\top}|^{2}H - \eta(H)\xi^{\top} + p\eta(H)\xi\} \\ &+& f_{3}\sum_{i}\{tr(P)s H + tr(P)t H - \tilde{g}(H, N e_{i})P e_{i} - \tilde{g}(H, N e_{i})N e_{i}\\ &+& 2\tilde{g}(e_{i}, s H)P e_{i} + 2\tilde{g}(e_{i}, s H)N e_{i}\}.
\end{eqnarray*}
It implies that
\begin{eqnarray*}
tr \big(R^{\star}(\cdot, H)\cdot\big) &=& - f_{1}pH 
+ f_{2}\{|\xi^{\top}|^{2}H - \eta(H)\xi^{\top} + p\eta(H)\xi\}+3f_3(PsH+NsH),
\end{eqnarray*}
by considering the anti-symmetry property of $\phi$, $tr(P)=0$ and $\tilde{g}(H, N e_{i}) = - \tilde{g}(t H, e_{i})$. 

Now, from value of $tr \big(R^{\star}(\cdot, H)\cdot\big) $ and equations \eqref{eqbiharmonic}, \eqref{fbiharmoniccondition}, we have result of the theorem by considering the tangential and normal parts.
\hfill$\square$\\ \\
Now, we have the following corollary if we assume different particular cases in Theorem \ref{thm2}.\linebreak
\begin{cor}\label{cor2}
Let $M^{p}$ be a submanifold of a generalized Sasakian space form $\widetilde{M}(f_1,f_2,f_3)$. 
\begin{enumerate}
\item If $M$ is invariant then $M$ is $f$-biharmonic if and only if
$$-\Delta^{\perp}H+\tr B(\cdot,A_H\cdot) + \frac{\Delta f}{f}H +  2\nabla^{\perp}_{\grad(\ln f)}H=pf_1H-f_2|\xi^{\top}|^2H-pf_2\eta(H)\xi^{\perp}$$
and
$$\frac{p}{2}\grad|H|^2+2\tr A_{\nabla^{\perp}H}(\cdot)-2A_{H}\grad(\ln f)=-2f_2(p-1)\eta(H)\xi^{\top}-6f_3PsH.$$

\item  If $M$ is anti-invariant then $M$ is $f$-biharmonic if and only if

\begin{eqnarray*}
-\Delta^{\perp}H+\tr B(\cdot,A_H\cdot) + \frac{\Delta f}{f}H +  2\nabla^{\perp}_{\grad(\ln f)}H=pf_1H\\-f_2|\xi^{\top}|^2H- pf_2\eta(H)\xi^{\perp}-3f_3NsH
\end{eqnarray*}
and
$$\frac{p}{2}\grad|H|^2+2\tr A_{\nabla^{\perp}H}(\cdot)-2A_{H} \grad(\ln f)=-2f_2(p-1)\eta(H)\xi^{\top}.$$
\item
If $\xi$ is normal to $M$ then $M$ is $f$-biharmonic if and only if
$$-\Delta^{\perp}H+\tr B(\cdot,A_H\cdot) + \frac{\Delta f}{f}H +  2\nabla^{\perp}_{\grad(\ln f)}H=pf_1H-pf_2\eta(H)\xi-3f_3NsH$$
and
$$ \frac{p}{2}\grad|H|^2+2\tr A_{\nabla^{\perp}H}(\cdot)-2A_{H} \grad(\ln f)=0.$$

\item If $\xi$ is tangent to $M$ then $M$ is $f$-biharmonic if and only if
$$-\Delta^{\perp}H+\tr B(\cdot,A_H\cdot) + \frac{\Delta f}{f}H +  2\nabla^{\perp}_{\grad(\ln f)}H=pf_1H-f_2H-3f_3NsH$$
and
$$ \frac{p}{2}\grad|H|^2+2\tr A_{\nabla^{\perp}H}(\cdot)-2A_{H} \grad(\ln f)=-6f_3PsH.$$
\item If $M$ is a hypersurface then $M$ is $f$-biharmonic if and only if
\begin{eqnarray*}
-\Delta^{\perp}H+\tr B(\cdot,A_H\cdot) + \frac{\Delta f}{f}H +  2\nabla^{\perp}_{\grad(\ln f)}H=(2nf_1 + 3f_3)H\\-f_2|\xi^{\top}|^2H-(2nf_2+3f_3)\eta(H)\xi^{\perp}
\end{eqnarray*}
and
$$ n\grad|H|^2+2\tr A_{\nabla^{\perp}H}(\cdot)-2A_{H} \grad(\ln f)=-(2(2n-1)f_1 +6f_3)\eta(H)\xi^{\top}.$$
\end{enumerate}
\end{cor}
{\bf Proof.} The proof is a direct consequence of Theorem \ref{thm3b} using the following facts.
\begin{enumerate}
\item If $M$ is invariant then $P=0$.
\item If $M$ is anti-invariant then $N=0$.
\item If $\xi$ is normal then $\eta(\grad f)=0$ and $M$ is anti-invariant which implies $P=0$.
\item If $\xi$ is tangent then $\eta(H)=0$.
\item If $M$ is a hypersurface then $sH=0$.
\end{enumerate}
Analogously to the case of generalized complex space forms (Proposition \ref{propB}), we can \linebreak obtain some curvature properties in some special cases by using characterizations of $f$-biharmonic submanifolds of generalized Sasakian space forms.
\begin{prop}\label{propscal}
$(1)$ Let $M^{2n}$ be a hypersurface of generalized Sasakian space form $\widetilde{M}(f_1,f_2,f_3)$ with non zero constant mean curvature $H$ and $\xi$ is tangent to $M$. Then 
 $M$ is proper $f$-biharmonic if and only if
\begin{center}
$|B|^2=2nf_1-f_2+3f_3-\frac{\Delta f}{f},\hspace{.2cm} and \hspace{.2cm}
A\,\grad f=0,$
\end{center}
or equivalently if and only if 
\begin{center}
$\scal_M=2n(2n-2)f_1+(4n-1)f_2-(2n-4)f_3+(2n-1)H^2+\frac{\Delta f}{f}H\hspace{.2cm}
and 
\hspace{.2cm}
A\,\grad f =0.$
\end{center}
$(2)$ There exists no proper $f$-biharmonic hypersurfaces with constant mean curvature and \linebreak constant scalar curvature so that $\xi$ is tangent.
\end{prop}
{\bf Proof.} Let $M$ be a $f$-biharmonic hypersurface of $\widetilde{M}(f_1,f_2,f_3)$ with non zero constant mean curvature and $\xi$ tangent to $M$. Then, from Corollary \ref{cor2}, we have
$$\left\{
\begin{array}{l}
-\Delta^{\perp}H+\tr B(\cdot,A_H\cdot) + \frac{\Delta f}{f}H +  2\nabla^{\perp}_{\grad (\ln f))}H\\ \\=(pf_1 + 3f_3)H-f_2|\xi^{\top}|^2H-(2nf_2+3f_3)\eta(H)\xi^{\perp},\\ \\
n\grad|H|^2+2\tr A_{\nabla^{\perp}H}(\cdot)-2A_{H} \grad (\ln f)=0.
\end{array}
\right.$$
Now, as per assumption, $\xi$ is tangent to $M$ which gives $\eta(H)=\eta(\nu)=0$. Therefore, we have
$$\phi^2\nu=-\nu+\eta(\nu)\xi=-\nu.$$
On the other hand, we have
\begin{eqnarray*}
\phi^2\nu&=&\phi(s\nu+t\nu)\\
&=&Ps\nu+Ns\nu+st\nu+t^2\nu.
\end{eqnarray*}
Hence, we get
\begin{equation}\label{minusnu}
-\nu=Ps\nu+Ns\nu+st\nu+t^2\nu.
\end{equation}
Moreover, since $\langle\phi\nu,\nu\rangle=\Omega(\nu,\nu)=0$, we have that $\phi\nu$ is tangent, i.e., $t\nu=0$. Thus, Equation \eqref{minusnu} becomes
$$-\nu=Ps\nu+Ns\nu,$$
and so $Ps=0$ and $Ns=-{\rm Id}$ by identification of tangential and normal parts. Using these results in the above $f$-biharmonic condition for the hypersurfaces of generalized Sasakian space forms, we have 
$$\left\{
\begin{array}{l}
\tr B(\cdot,A_H\cdot)  =(2nf_1 + 3f_3)H-f_2|\xi^{\top}|^2H - \frac{\Delta f}{f}H,\\ \\
A_{H} \grad(\ln f)=0.
\end{array}
\right.$$
Hence, the second equation is trivial and the first becomes
$$\tr B(\cdot,A_H\cdot)=2nf_1H-f_2H+3f_3H-\frac{\Delta f}{f}H,$$
or equivalently
$$|B|^2=2nf_1-f_2+3f_3-\frac{\Delta f}{f},$$
since $\tr B(\cdot,A_H\cdot)=|B|^2H$ and $H$ is a non zero constant.\\
\\
Similarly, using Gauss formula for second part, we have
\begin{eqnarray*}
Scal_{M} &=& \sum_{i, j}\tilde{g}(R^{\star}(e_{i}, e_{j})e_{j}, e_{i}) - |B|^{2} - pH^{2}\\
&=& \sum_{i, j}f_{1}\{\tilde{g}(e_{j}, e_{j})\tilde{g}(e_{i}, e_{i}) - \tilde{g}(e_{i}, e_{j})\tilde{g}(e_{j}, e_{i})\}
+ \sum_{i, j}f_{2}\{\eta(e_{i})\eta(e_{j})\tilde{g}(e_{j}, e_{i})\\ &-& \eta(e_{j})\eta(e_{j})\tilde{g}(e_{i}, e_{i}) + \tilde{g}(e_{i}, e_{j})\eta(e_{j})\tilde{g}(\xi, e_{i}) - \tilde{g}(e_{j}, e_{j})\eta(e_{i})\tilde{g}(\xi, e_{i})\}\\
&+& \sum_{i, j}f_{3}\{\tilde{g}(e_{i}, \phi e_{j})\tilde{g}(\phi e_{j}, e_{i}) - \tilde{g}(e_{j}, \phi e_{j})\tilde{g}(\phi e_{i}, e_{i}) + 2\tilde{g}(e_{i}, \phi e_{j})\tilde{g}(\phi e_{j}, e_{i})\}\\ &-& |B|^{2} - pH^{2}
= 2n(2n - 1)f_{1} + 2(2n - 1)f_{2} - (2n - 1)f_{3} - |B|^{2} - pH^{2}.
\end{eqnarray*}
Using the value of $|B|^{2}$ obtain in the first part of the proof, we get the required result, that is,
$$\scal_M=2n(2n-2)f_1+(4n-1)f_2-(2n-4)f_3+(2n-1)H^2+\frac{\Delta f}{f}H.$$
Moreover, since $\grad (\ln f)=\frac{1}{f}\grad f$  and $A_H=HA$ with $H$ is a positive constant, the equation $A_{H} \grad(\ln f)=0$ reduces to $A\,\grad f=0$. This concludes the proof.
\hfill$\square$\\ \\
Now, from this proposition, we can prove the following non-existence result.
\begin{cor}
Let $M^{2n}$ be a constant mean curvature hypersurface of generalized Sasakian space form $\widetilde{M}(f_1,f_2,f_3)$ with $\xi$ tangent. If the functions $f_1, f_2, f_3$ satisfy the inequality \linebreak $2nf_1-f_2+3f_3\leqslant \frac{(\Delta f)}{f}$ on $M$ then $M$ is not biharmonic.\\ 
 In particular, there exists no proper $f$-biharmonic CMC hypersurface with $\xi$ tangent and $f$ \linebreak satisfying 
\begin{itemize}
\item   \hspace{.2cm} $\tilde{c}\leqslant\frac{4}{2n+2}[\frac{\Delta f}{f}-\frac{6n-2}{4}]$ in a Sasakian space form $\widetilde{M}^{2n+1}_S(\tilde{c})$.
\item  \hspace{.2cm}$\tilde{c}\leqslant\frac{4}{2n+2}[\frac{\Delta f}{f}+\frac{6n-2}{4}]$ in a Kenmotsu space form $\widetilde{M}^{2n+1}_K(\tilde{c})$.
\item   \hspace{.2cm} $\tilde{c}\leqslant\frac{4}{2n+2}\frac{\Delta f}{f}$ in a cosymplectic space form $\widetilde{M}^{2n+1}_C(\tilde{c})$.
\end{itemize}
\end{cor}
{\bf Proof:} As per assumption, $M$ is a hypersurface of $\widetilde{M}(f_1,f_2,f_3)$ with non zero constant mean curvature $H$ and $\xi$ tangent to $M$. From Proposition \ref{propscal}, $M$ is $f$-biharmonic if and only if its second fundamental form $B$ satisfies
$|B|^2=2nf_1-f_2+3f_3-\frac{\Delta f}{f}.$
In other words, this is not possible if 
\begin{equation}\label{condf1f2f3}
2nf_1-f_2+3f_3\leqslant\frac{\Delta f}{f}.
\end{equation}\\
Now, $f_1=\frac{\tilde{c}+3}{4}$ and $f_2=f_3=\frac{\tilde{c}-1}{4}$ if $\widetilde{M}(f_1,f_2,f_3)$ is a Sasakian space form where $\tilde{c}$ is $\phi$-sectional curvature. Therefore, the inequality $2nf_1-f_2+3f_3\leqslant\frac{\Delta f}{f}$ reduces to $\tilde{c}\leqslant\frac{4}{2n+2}[\frac{\Delta f}{f}-\frac{6n-2}{4}]$. Similarly, we have $f_1=\frac{\tilde{c}-3}{4}$ and $f_2=f_3=\frac{\tilde{c}+1}{4}$ ({\it resp.} $f_1=f_2=f_3=\frac{\tilde{c}}{4}$) for the Kenmotsu ({\it resp. } cosymplectic) case and the inequality $2nf_1-f_2+3f_3\leqslant\frac{\Delta f}{f}$ reduces to $\tilde{c}\leqslant\frac{4}{2n+2}[\frac{\Delta f}{f}+\frac{6n-2}{4}]$ ({\it resp.} $\tilde{c}\leqslant\frac{4}{2n+2}\frac{\Delta f}{f}$).
\hfill $\square$\\ \\
Now, we have the following proposition analogous to complex case.
\begin{thm}\label{thmSKC}
Let $M^{q}$  be a submanifold of Sasakian (Kenmotsu or cosymplectic) space form $\widetilde{M}^{2n+1}_S(\tilde{c})$ ({\it resp.} $\widetilde{M}^{2n+1}_K(\tilde{c})$ or $\widetilde{M}^{p+1}_C(\tilde{c})$) with constant mean curvature $H$ so that $\xi$ and $\phi H$ are tangent.  Further, we consider $F(f,q,\tilde{c})$ the function defined on $M$ by
$$ F(f,q,\tilde{c}) = qf_1-f_2+3f_3-\frac{\Delta f}{f}=\left\{\begin{array}{ll}
\frac{(q+2)\tilde{c}}{4}+\frac{(3q-2)}{4}-\frac{\Delta f}{f}&\text{for}\hspace{.2cm}\widetilde{M}^{p+1}_S(\tilde{c}),\\ \\
\frac{(q+2)\tilde{c}}{4}-\frac{(3q-2)}{4}-\frac{\Delta f}{f}&\text{for}\hspace{.2cm} \widetilde{M}^{p+1}_K(\tilde{c}),\\ \\
\frac{(q+2)\tilde{c}}{4}-\frac{\Delta f}{f}&\text{for}\hspace{.2cm} \widetilde{M}^{p+1}_C(\tilde{c}).
\end{array}
\right.$$
Then we have the following observations.
\begin{enumerate} 
\item If \hspace{.2 cm}$\displaystyle\inf_M F(f,q,\tilde{c})$ is non-positive then $M$ is not $f$-biharmonic. \\
\item
If \hspace{.2 cm}$\displaystyle\inf_MF(f,q,\tilde{c})$ is positive and $M$ is proper $f$-biharmonic then 
$$0<|H|^2 \leqslant \frac1q\inf_{M}F(f,q,\tilde{c}).$$

\end{enumerate}\end{thm}
{\bf Proof:} As $M$ is proper $f$-biharmonic submanifold with constant mean curvature $H$ and $\xi$ tangent to $M$, so we get form Corollary \ref{cor2} that
$$\left\{
\begin{array}{l}
-\Delta^{\perp}H+\tr B(\cdot,A_H\cdot) + \frac{2}{f}\nabla^{\perp}_{\grad f}H+ \frac{\Delta f}{f}H=qf_1H-f_2H-3f_3NtH,\\ \\
2\tr A_{\nabla^{\perp}H}(\cdot)-2A_{H} \grad(\ln f)=-6f_3PtH.
\end{array}
\right.$$
Now, considering $\phi H$ is tangent implies that $sH=0$. Again applying $\phi$ gives that $\phi^2H=PtH+NtH$. But from $\phi^2H=-H+\eta(H)\xi$ and $\xi$ is tangent, we have $\phi^2H=-H$. Therefore, comparing tangential and normal parts, we get $PtH=0$ and $NtH=-H$. Using these facts in the above equation, we get
$$\left\{
\begin{array}{l}
-\Delta^{\perp}H+\tr B(\cdot,A_H\cdot)=qf_1H-f_2H+3f_3H -\frac{\Delta f}{f}H,\\ \\
2\tr A_{\nabla^{\perp}H}(\cdot)-2A_{H} \grad(\ln f)=0.
\end{array}
\right.$$
Now, considering $\nu$ as an real eigenvalue of the eigenfunction $f$ corresponding to Laplacian operator $\Delta$, i.e., $\frac{\Delta f}{f} = \nu$,
from first equation, we have

\begin{eqnarray*}-\Delta^{\perp}H+\tr B(\cdot,A_H\cdot)&=&qf_1H-f_2H+3f_3H-\nu H\\
&=&F(f,q,\tilde{c})H.
\end{eqnarray*}
Taking scalar product by $H$, we get
 $$-\left\langle\Delta^{\perp}H,H\right\rangle+\left\langle\tr B(\cdot,A_H\cdot),H\right\rangle=F(f,q,\tilde{c})|H|^2.$$
Using the facts $\left\langle\tr B(\cdot,A_H\cdot),H\right\rangle=|A_H|^2$, $|H|$ is a constant and the B\"ochner formula, i.e., $\frac{1}{2}\Delta|H|^2=\left\langle\Delta^{\perp}H,H\right\rangle-|\nabla^{\perp}H|^2$ in the above equation, we have
$$|A_H|^2+|\nabla^{\perp}H|^2=F(f,q,\tilde{c})|H|^2.$$ 
Now, this equation reduces to
$$F(f,q,\tilde{c})|H|^2=|A_H|^2+|\nabla^{\perp}H|^4\geqslant q|H|^2+|\nabla^{\perp}H|^2\geqslant q|H|^4,$$
by considering the Cauchy-Schwarz inequality $|A_H|^2\geqslant \frac{1}{q}\tr(A_H)=q|H|^4$. It implies that
$$F(f,q,\tilde{c})\geqslant q|H|^2,$$
as $|H|$ is a positive constant. This proves the two assertions of the theorem.
\hfill$\square$\\ \\
Now, we have the analogous result replacing the assumption that $\phi H$ is tangent by $\phi H$ is normal. Namely, we have:
\begin{prop}
Let $\psi: M^{q}\rightarrow\widetilde{M}^{p+1}_S(\tilde{c})$ ({\it resp.} $\widetilde{M}^{p+1}_K(\tilde{c})$ or $\widetilde{M}^{p+1}_C(\tilde{c})$) be a submanifold of Sasakian (Kenmotsu or cosymplectic) space form with constant mean curvature $H$ so that $\xi$ is tangent and $\phi H$ is normal. Further, we consider $F(f,q,\tilde{c})$ the function defined on $M$ by
$$ G(f,q,\tilde{c}) = qf_1-f_2-\frac{\Delta f}{f}=\left\{\begin{array}{ll}
\frac{(q-1)\widetilde{c}}{4}+\frac{(3q+1)}{4}-\frac{\Delta f}{f}&\text{for}\hspace{.2cm}\widetilde{M}^{p+1}_S(\tilde{c}),\\ \\
\frac{(q-1)\widetilde{c}}{4}-\frac{(3q+1)}{4}-\frac{\Delta f}{f}&\text{for}\hspace{.2cm} \widetilde{M}^{p+1}_K(\tilde{c}),\\ \\
\frac{(q-1)\widetilde{c}}{4}-\frac{\Delta f}{f}&\text{for}\hspace{.2cm} \widetilde{M}^{p+1}_C(\tilde{c}).
\end{array}
\right.$$
Then we have the following observations.
\begin{enumerate} 
\item If $\displaystyle\inf_M G(f,q,\tilde{c})$ is non-positive then $M$ is not $f$-biharmonic. \\
\item
If $\displaystyle\inf_MG(f,q,\tilde{c})$ is positive and $M$ is proper $f$-biharmonic then 
$$0<|H|^2 \leqslant \frac1q\inf_{M}G(f,q,\tilde{c}).$$

\end{enumerate}

\end{prop}
{\bf Proof:} 
Now, in this case, $M$ is proper $f$-biharmonic submanifold with $\xi$ is tangent and $\phi H$ is normal. Normality of $\phi H$ implies that $sH=0$. Therefore, from Corollary \ref{cor2}, we have
\begin{eqnarray*}-\Delta^{\perp}H+\tr B(\cdot,A_H\cdot)&=&qf_1H-f_2H\\
&=&G(f,q,\tilde{c})H.
\end{eqnarray*}
Similarly, as in the previous theorem, taking the scalar product by $H$ and using the B\"ochner formula and then with the help of the Cauchy-Schwarz inequality, we get
$$G(f,q,\tilde{c})|H|^2=|A_H|^2+|\nabla^{\perp}H|^4\geqslant q|H|^2+|\nabla^{\perp}H|^2\geqslant q|H|^4.$$
It easily provides the inequality $G(f,q,\tilde{c})\geqslant q|H|^2$, since $|H|$ is a positive constant. We get $0<|H|^2 \leqslant \frac1q\inf_{M}G(f,q,\tilde{c}),$ which concludes the proof.\hfill$\square$

\section{\textbf {Bi-$f$-harmonic submanifolds}}\label{sec4}
In this section, we consider bi-$f$-harmonic submanifolds, which are, as we mention in the \linebreak introduction, different from the $f$-biharmonic submanifolds studied above. 
\subsection{A general necessary and sufficient condition}
We begin by giving this general result which gives the necessary and sufficient condition for a submanifold to be bi-$f$-harmonic
\begin{thm}\label{thmgene}
Let $(M^n,g)$ be a Riemannian manifold isometrically immersed into another Riemannian manifold $(N,h)$. Let $f$ be a smooth function on $M$. Then, $M$ is a bi-$f$-harmonic submanifold of $N$ if and only if the following two equations hold:
\begin{enumerate}
\item
$$ \begin{array}{l}
nf^2\Delta^{\perp}H+nf^2\tr B(\cdot,A_H\cdot)-nf(\Delta f)H-3n\nabla^{\perp}_{\grad f}H\\
-f\tr B(\cdot,\nabla_{\cdot}\grad f)-f\tr\nabla_{\cdot}B(\cdot,\grad f)-n|\grad f|^2H-B(\grad f,\grad f)\\
=-nf^2\tr\left(\overline{R}(\cdot,H)\cdot\right)^{\perp}-f\tr\left(\overline{R}(\cdot,\grad f)\cdot\right)^{\perp},
\end{array}$$
\item $$\begin{array}{l}
\dfrac{n^2f^2}{2}\grad|H|^2+2n^2f^2\tr(A_{\nabla^{\perp}_{\cdot}H\cdot})+3nfA_H\grad f\\
+fRic_M(\grad f)+f\grad(\Delta f)+f\tr(A_{B(\cdot,\grad f)}(\cdot))-\dfrac12\grad(|\grad f|^2)\\
=-2nf^2\tr\left(\overline{R}(\cdot,H)\cdot\right)^{\top}-f\tr\left(\overline{R}(\cdot,\grad f)\cdot\right)^{\top}.
\end{array}$$
\end{enumerate}

\end{thm}
\textbf{Proof:} We recall that $M$ is a bi-$f$-harmonic submanifold if and only if 
$$\tau^{2}_f(\psi)=fJ^{\psi}(\tau_f(\psi)) - \nabla^{\psi}_{grad f}\tau_f(\psi) = 0,$$
where $J^{\psi}$ is defined by $$J^{\psi}(X) = -[Tr_{g}\nabla^{\psi}\nabla^{\psi}X - \nabla^{\psi}_{\nabla^{M}}X - R^{N}(d\psi, X)d\psi]$$
and $\tau_f(\psi)=f\tau(\psi)+d\psi(\grad f)$. Since we are in the case of submanifolds, for a sake of \linebreak compactness, we will omit the map $\psi$ and we will denote $\nabla^{\psi}$ by $\overline{\nabla}$ as the Levi-Civita \linebreak connection on $N$. Hence, we have 
 $$\tau_f(\psi)=f\tau(\psi)+d\psi(\grad f)=nfH+\grad f.$$
 Taking $\{e_1,\cdots,e_n\}$ a normal frame of $T_pM$ for a fixed point $p\in M$, we get
 \begin{eqnarray}\label{trR1}
 Tr_{g}\Big(\nabla^{\psi}\nabla^{\psi}\tau_f(\psi) - \nabla^{\psi}_{\nabla^{M}}\tau_f(\psi)\Big)&=&\sum_{i=1}^n\nabla^{\psi}_{e_i}\nabla^{\psi}_{e_1}\tau_f(\psi)- \nabla^{\psi}_{\nabla^{M}_{e_1}e_1}\tau_f(\psi)\nonumber\\
 &=&\sum_{i=1}^n\overline{\nabla}_{e_i}\overline{\nabla}_{e_1}(nfH+\grad f).
 \end{eqnarray}
 First, we compute 
 \begin{eqnarray}\label{trR2}
 \sum_{i=1}^n\overline{\nabla}_{e_i}\overline{\nabla}_{e_1}(nfH)&=&n\sum_{i=1}^n\overline{\nabla}_{e_1}(e_i(f)H+f\overline{\nabla}_{e_1}H)\nonumber\\
 &=&n\sum_{i=1}^n\Big(  e_i(e_i(f))H+2e_i(f)\overline{\nabla}_{e_1}H+f\overline{\nabla}_{e_1}\overline{\nabla}_{e_1}H\Big)\nonumber\\
 &=&-n\Delta f+2\overline{\nabla}_{\grad f}H+nf\sum_{i=1}^n\overline{\nabla}_{e_1}\overline{\nabla}_{e_1}H.
 \end{eqnarray}
 Now, we give this first lemma.
 \begin{lemma}\label{lemgene1}
 We have 
 \begin{eqnarray*}
 \sum_{i=1}^n\overline{\nabla}_{e_1}\overline{\nabla}_{e_1}H=-\frac{n}{2}\grad |H|^2-\tr B(\cdot,A_H\cdot)-2\tr(A_{\nabla^{\perp}_{\cdot}H}\cdot)-\Delta^{\perp}H+\tr(\overline{R}(\cdot,H)\cdot)^{\top}.
 \end{eqnarray*}
 \end{lemma}
 \textbf{Proof:}
 We have
 \begin{eqnarray*}
 \overline{\nabla}_{e_1}\overline{\nabla}_{e_1}H&=& \overline{\nabla}_{e_1}(\nabla^{\perp}_{e_1}H-A_He_i)\\
 &=&\nabla^{\perp}_{e_1}\nabla^{\perp}_{e_i}H-A_{\nabla^{\perp}_{e_i}H}-\nabla_{e_i}(A_He_i)-B(A_He_i,e_i)
 \end{eqnarray*}
 Hence, summing over $i$, we get
 \begin{equation}\label{lem2eq1}
 \sum_{i=1}^n \overline{\nabla}_{e_1}\overline{\nabla}_{e_1}H = -\Delta^{\perp}H-\tr(A_{\nabla^{\perp}_{(\cdot)}H}(\cdot))-\tr(\nabla_{(\cdot)}A_H(\cdot))-\tr B(\cdot,A_H(\cdot)).
 \end{equation}
 Moreover, we have
 \begin{eqnarray}\label{lem2eq2}
\tr(\nabla_{(\cdot)}A_H(\cdot))&=&\sum_{i=1}^n\nabla_{e_1}(A_He_i)
 =\sum_{i,j=1}^ng(\nabla_{e_1}(A_He_i),e_j)e_j\nonumber\\
 &=&\sum_{i,j=1}^ne_ig(A_He_i,e_j)e_j\nonumber
 =\sum_{i,j=1}^ne_ig(B(e_i,e_j),H)e_j\nonumber\\
 &=&\sum_{i,j=1}^ne_ig(\overline{\nabla}_{e_j}e_i,H)e_j\nonumber\\
 &=&\sum_{i,j=1}^n\Big(g(\overline{\nabla}_{e_i}\overline{\nabla}_{e_j}e_i,H)e_j+ g(\overline{\nabla}_{e_j}e_i,\overline{\nabla}_{e_i}H)e_j\Big)\nonumber\\
 &=&\sum_{i,j=1}^n\Big(g(\overline{\nabla}_{e_i}\overline{\nabla}_{e_j}e_i,H)e_j+ g(B(e_j,e_i),\overline{\nabla}_{e_i}H)e_j\Big)\nonumber\\
 &=&\sum_{i,j=1}^ng(\overline{\nabla}_{e_i}\overline{\nabla}_{e_j}e_i,H)e_j+ \sum_{i=1}^nA_{\nabla^{\perp}_{e_i}H}e_i\nonumber\\
 &=&\sum_{i,j=1}^ng(\overline{\nabla}_{e_i}\overline{\nabla}_{e_j}e_i,H)e_j+\tr(A_{\nabla^{\perp}_{\cdot}H}\cdot).
 \end{eqnarray}
 Now, we have
 \begin{eqnarray*}
 \sum_{i,j=1}^ng(\overline{\nabla}_{e_i}\overline{\nabla}_{e_j}e_i,H)e_j&=&\sum_{i,j=1}^ng\Big(\overline{R}(e_i,e_j)e_i+\overline{\nabla}_{e_j}\overline{\nabla}_{e_i}e_i+\overline{\nabla}_{[e_i,e_j]}e_i,H\Big)e_j.
 \end{eqnarray*}
 Since, the frame $\{e_1,\cdot\cdot\cdot, e_n\}$ is normal, we have $[e_i,e_j]=0$. Moreover, we have
 $$\sum_{i=1}^n\overline{\nabla}_{e_j}\overline{\nabla}_{e_i}e_i=n\overline{\nabla}_{e_j}H.$$
 Hence, we get
 \begin{eqnarray}\label{lem2eq3}
 \sum_{i,j=1}^ng(\overline{\nabla}_{e_i}\overline{\nabla}_{e_j}e_i,H)e_j&=&-\tr(\overline{R}(\cdot,H)\cdot)^{\top}+ n\sum_{j=1}^ng(\overline{\nabla}_{e_j}H,H)e_j\nonumber\\
 &=&-\tr(\overline{R}(\cdot,H)\cdot)^{\top}+\frac{n}{2}\grad|H|^2.
 \end{eqnarray}
 Putting \eqref{lem2eq3} into \eqref{lem2eq2}, we get
 \begin{equation}
 \tr(\nabla_{(\cdot)}A_H(\cdot))=\tr(A_{\nabla^{\perp}_{\cdot}H}\cdot)-\tr(\overline{R}(\cdot,H)\cdot)^{\top}+\frac{n}{2}\grad|H|^2,
 \end{equation}
 and finally, reporting this in \eqref{lem2eq1}, we get
  \begin{eqnarray*}
 \sum_{i=1}^n\overline{\nabla}_{e_1}\overline{\nabla}_{e_1}H=-\frac{n}{2}\grad|H|^2-\tr B(\cdot,A_H\cdot)-2\tr(A_{\nabla^{\perp}_{\cdot}H}\cdot)-\Delta^{\perp}H+\tr(\overline{R}(\cdot,H)\cdot)^{\top},
 \end{eqnarray*}
 which concludes the proof of the lemma.
 \hfill$\square$\\
 We now state this second lemma.
  \begin{lemma}\label{lemgene2}
 We have 
 \begin{eqnarray*}
 \sum_{i=1}^n\overline{\nabla}_{e_i}\overline{\nabla}_{e_i}\grad f&=&\grad(\Delta f)+2{\rm Ric}_M(\grad f)-\tr(\overline{R}(\cdot,\grad f)\cdot)\\
 &&+\tr(B(\cdot,\nabla_{\cdot}\grad f))+\tr(\nabla^{\perp}_{\cdot}B(\cdot,\grad f))-\tr(A_{B(\cdot,\grad f)}(\cdot)).
 \end{eqnarray*}
 \end{lemma}
 \textbf{Proof:} We have
 \begin{eqnarray*}
 \overline{\nabla}_{e_i}\overline{\nabla}_{e_i}\grad f&=& \overline{\nabla}_{e_i}\big(\nabla_{e_i}\grad f+B(e_i,\grad f)\big)\\
 &=&\nabla_{e_i}\nabla_{e_i}\grad f+B(e_i,\nabla_{e_i}\grad f)+\nabla^{\perp}_{e_i}B(e_i,\grad f)-A_{B(e_i,\grad f)}(e_i)\\
 &=&\nabla_{e_i}\nabla_{e_i}\grad f+\tr(B(\cdot,\nabla_{\cdot}\grad f))+\tr(\nabla^{\perp}_{\cdot}B(\cdot,\grad f))-\tr(A_{B(\cdot,\grad f)}(\cdot)).
 \end{eqnarray*}
Moreover, we conclude by using the classical fact that (see \cite{} for instance)
 $$\sum_{i=1}^n\nabla_{e_i}\nabla_{e_i}\grad f=\grad(\Delta f)+2{\rm Ric}_M(\grad f)-\tr(\overline{R}(\cdot,\grad f)\cdot).$$
\hfill$\square$\\
Finally, we have this last elementary lemma.
\begin{lemma}\label{lemgene3}
We have
\begin{eqnarray*}\overline{\nabla}_{grad (f)}\tau_f(\psi)&=&n|\grad f|^2H-nfA_H(\grad f)+nf\nabla^{\perp}_{\grad f}H\\
&&+\frac12\grad(|\grad f|^2)+B(\grad f,\grad f).
\end{eqnarray*}
\end{lemma}
\textbf{Proof:} We have
\begin{eqnarray*}\overline{\nabla}_{grad f}\tau_f(\psi)&=&\overline{\nabla}_{grad f}\big(nfH+\grad f\big)\\
&=&n|\grad f|^2H+nf\nabla^{\perp}_{\grad f}H-nfA_H(\grad f)\\
&&+\nabla_{\grad f}\grad f+B(\grad f,\grad f).
\end{eqnarray*}
Using the fact that $\nabla_{\grad f}\grad f=\frac12\grad (|\grad f|^2)$, we get the desired identity.
\hfill$\square$\\\\
Now, we can finish the proof of Theorem \ref{thmgene}. Recall that $M$ is a bi-$f$-harmonic submanifold\linebreak if and only if 
\begin{eqnarray*}
\tau_f^2(\phi)=fJ^{\psi}(\tau_f(\psi)) - \nabla^{\psi}_{grad (f)}\tau_f(\psi)=0.
\end{eqnarray*}
From \eqref{trR1} and \eqref{trR2}, we have
\begin{eqnarray*}
\tau_f^2(\phi)=-nf\Delta f+2f\overline{\nabla}_{\grad f}H+nf^2\sum_{i=1}^n\overline{\nabla}_{e_1}\overline{\nabla}_{e_1}H+f\sum_{i=1}^n\overline{\nabla}_{e_1}\overline{\nabla}_{e_1}\grad f- \nabla^{\psi}_{grad f}\tau_f(\psi) .
\end{eqnarray*}
Replacing the last three terms in the right-hand side using, respectively, Lemmas \ref{lemgene1}, \ref{lemgene2} and \ref{lemgene3}, we obtain
\begin{eqnarray*}
\tau_f^2(\phi)&=&-nf\Delta f+2f\overline{\nabla}_{\grad f}H-nf^2\frac{n}{2}\grad|H|^2-nf^2\tr B(\cdot,A_H\cdot)-2nf^2\tr(A_{\nabla^{\perp}_{\cdot}H}\cdot)\\&&-nf^2\Delta^{\perp}H+nf^2\tr(\overline{R}(\cdot,H)\cdot)^{\top}+f\grad(\Delta f)+2f{\rm Ric}_M(\grad f)\\&&-f\tr(\overline{R}(\cdot,\grad f)\cdot)+f\tr(B(\cdot,\nabla_{\cdot}\grad f))+f\tr(\nabla^{\perp}_{\cdot}B(\cdot,\grad f))\\&&-f\tr(A_{B(\cdot,\grad f)}(\cdot))-n|\grad f|^2H+nfA_H(\grad f)-nf\nabla^{\perp}_{\grad f}H\\
&&-\frac12\grad(|\grad f|^2)-B(\grad f,\grad f).
\end{eqnarray*}
Finally decomposing the fact that $\tau_f^2(\phi)=0$ into tangent and normal parts, we get the two identities of the theorem. This conludes the proof of Theorem \ref{thmgene}.\hfill$\square$
 \subsection{Bi-$f$-harmonic submanifolds of generalized complex space forms}
 In this section, using the general bi-$f$-harmonicity condition of Theorem \ref{thmgene}, we give the necessary and sufficient condition for submanifold of generalized complex space forms to be bi-$f$-harmonic. Namely, we have the following theorem.
\begin{thm}\label{thm1b}
Let $M^{n}$, $n<4$, be a submanifold of generalized complex space form $N(\alpha,\beta)$ with second fundamental form $B$, shape operator $A$, mean curvature $H$ and a positive $C^{\infty}$-differentiable function $f$ on $M$. Then $M$ is bi-$f$-harmonic submanifold of $N(\alpha,\beta)$ if and only the following two equations are satisfied
\begin{enumerate}
\item
$$ \begin{array}{l}
nf^2\Delta^{\perp}H+nf^2\tr B(\cdot,A_H\cdot)-nf(\Delta f)H-3n\nabla^{\perp}_{\grad f}H\\
-f\tr B(\cdot,\nabla_{\cdot}\grad f)-f\tr\nabla_{\cdot}B(\cdot,\grad f)-n|\grad f|^2H-B(\grad f,\grad f)\\
=n^2f^2\alpha H-3nf^2\beta klH-3f\beta kj\grad f,
\end{array}$$
\item $$\begin{array}{l}
\dfrac{n^2f^2}{2}\grad|H|^2+2n^2f^2\tr(A_{\nabla^{\perp}_{\cdot}H\cdot})+3nfA_H\grad f\\
+fRic_M(\grad f)+f\grad(\Delta f)+f\tr(A_{B(\cdot,\grad f)}(\cdot))-\dfrac12\grad(|\grad f|^2)\\
=-6nf^2\beta jlH+2f(n-1)\alpha\grad f-6f\beta j^2\grad f.
\end{array}$$
\end{enumerate}

\end{thm}
\noindent
{\bf Proof:} This theorem is a direct consequence of Theorem \ref{thmgene} taking into account that the\linebreak curvature tensor of the generalized complex space form $N(\alpha,\beta)$ is given by $R=\alpha R_1+\beta R_2$, with $R_1$ and $R_2$ defined in Section \ref{sec2}. First, we have from \eqref{curvfbiharm}
$$\tr\left(R(\cdot,H)\cdot\right)=-n\alpha H + \beta(3jlH + 3klH).$$
Moreover, we need to compute $\tr\left(\overline{R}(\cdot,\grad f)\cdot\right)$. We have the following lemma
\begin{lemma}
We have
$$\tr\left(\overline{R}(\cdot,\grad f)\cdot\right)^{\top}=-(n-1)\alpha\grad f+3\beta j^2\grad f$$
and 
$$\tr\left(\overline{R}(\cdot,\grad f)\cdot\right)^{\perp}=3\beta kj\grad f.$$
\end{lemma}
{\bf Proof:} We have
\begin{eqnarray*}
\tr\left(\overline{R}(\cdot,\grad f)\cdot\right)&=&\alpha \sum_{i=1}^n\left(\left\langle\grad f,e_i\right\rangle e_i-\left\langle e_i,e_i\right\rangle\grad f\right)\\
&&+\beta\sum_{i=1}^n\left( \left\langle J\grad f,e_i\right\rangle Je_i -\left\langle Je_i,e_i\right\rangle J\grad f+2\left\langle J\grad f,e_i\right\rangle Je_i\right)\\
&=&-(n-1)\alpha\grad f+3\beta J(J\grad f)^{\top}\\
&=&-(n-1)\alpha\grad f+3\beta j^2\grad f+kj\grad f.
\end{eqnarray*}
We conclude the proof of the lemma by identifying tangential and normal parts.
\hfill$\square$\\ \\
Using this lemma and reporting into Theorem \ref{thmgene}, we get the desired identities. This concludes the proof of Theorem \ref{thm1b}.
\hfill$\square$
\begin{thm}\label{thm2b}
Let $\psi: M^{n}\rightarrow \mathbb{M}^N_{\mathbb{C}}(c)$,  $n\leqslant 2N$, be a submanifold of the complex space form $\mathbb{M}^N_{\mathbb{C}}(c)$ of complex dimension $N$ and constant holomorphic sectional curvature $c$, with second fundamental form $B$, shape operator $A$, mean curvature $H$ and a positive $C^{\infty}$-differentiable function $f$ on $M$. Then $M$ is bi-$f$-harmonic submanifold of $\mathbb{M}^N_{\mathbb{C}}(c)$ if and only the following two equations are satisfied
\begin{enumerate}
\item
$$ \begin{array}{l}
nf^2\Delta^{\perp}H+nf^2\tr B(\cdot,A_H\cdot)-nf(\Delta f)H-3n\nabla^{\perp}_{\grad f}H\\
-f\tr B(\cdot,\nabla_{\cdot}\grad f)-f\tr\nabla_{\cdot}B(\cdot,\grad f)-n|\grad f|^2H-B(\grad f,\grad f)\\
=n^2f^2c H-3nf^2c klH-3fc kj\grad f,
\end{array}$$
\item $$\begin{array}{l}
\dfrac{n^2f^2}{2}\grad|H|^2+2n^2f^2\tr(A_{\nabla^{\perp}_{\cdot}H\cdot})+3nfA_H\grad f\\
+fRic_M(\grad f)+f\grad(\Delta f)+f\tr(A_{B(\cdot,\grad f)}(\cdot))-\dfrac12\grad(|\grad f|^2)\\
=-6nf^2c jlH+2f(n-1)c\grad f-6fc j^2\grad f.
\end{array}$$
\end{enumerate}
\end{thm}
\noindent
{\bf Proof:} The proof is similar to the one of Theorem \ref{thm1b} with the only difference that $\alpha=\beta=c$ and $n<2N$ instead of $n<4$.
\hfill$\square$
\\
Now, we consider some particular cases where these conditons become simpler. Namely, we have:

\begin{cor}\label{corlag}
Let $\psi: M^{p}\rightarrow N(\alpha,\beta)$, $p<4$ be a submanifold of generalized complex space form $N(\alpha,\beta)$ with second fundamental form $B$, shape operator $A$, mean curvature $H$ and a positive $C^{\infty}$-differentiable function $f$ on $M$.
\begin{enumerate}

\item If $M$ is a hypersurface of $N(\alpha,\beta)$ with constant mean curvature, then $M$ is bi-$f$-harmonic if and only if

$$ \left\{\begin{array}{l}

-f\tr B(\cdot,\nabla_{\cdot}\grad f)-f\tr\nabla_{\cdot}B(\cdot,\grad f)-n|\grad f|^2H-B(\grad f,\grad f)\\
=n^2f^2\alpha H+3nf^2\beta H-nf^2H|B|^2+nf(\Delta f)H,\\ \\
fRic_M(\grad f)+f\grad(\Delta f)+f\tr(A_{B(\cdot,\grad f)}(\cdot))-\dfrac12\grad(|\grad f|^2)\\
=2f(n-1)\alpha\grad f-6f\beta j^2\grad f-3nfA_H\grad f.
\end{array}\right.$$

\item If $M$ is a complex surface of $N(\alpha,\beta)$ with parallel mean curvature, then $M$ is bi-$f$-harmonic if and only if

$$\left\{ \begin{array}{l}

-f\tr B(\cdot,\nabla_{\cdot}\grad f)-f\tr\nabla_{\cdot}B(\cdot,\grad f)-n|\grad f|^2H-B(\grad f,\grad f)\\
=n^2f^2\alpha H-nf^2\tr B(\cdot,A_H\cdot)+nf(\Delta f)H,\\ \\

+fRic_M(\grad f)+f\grad(\Delta f)+f\tr(A_{B(\cdot,\grad f)}(\cdot))-\dfrac12\grad(|\grad f|^2)\\
=2f(n-1)\alpha\grad f-6f\beta j^2\grad f-3nfA_H\grad f.
\end{array}\right.$$

\item If $M$ is a Lagrangian surface of $N(\alpha,\beta)$ with parallel mean curvature, then $M$ is bi-$f$-harmonic if and only if

$$ \left\{\begin{array}{l}

-f\tr B(\cdot,\nabla_{\cdot}\grad f)-f\tr\nabla_{\cdot}B(\cdot,\grad f)-n|\grad f|^2H-B(\grad f,\grad f)\\
=n^2f^2\alpha H+3nf^2\beta H-nf^2\tr B(\cdot,A_H\cdot)+nf(\Delta f)H\\\\ 

+fRic_M(\grad f)+f\grad(\Delta f)+f\tr(A_{B(\cdot,\grad f)}(\cdot))-\dfrac12\grad(|\grad f|^2)\\
=2f(n-1)\alpha\grad f-3nfA_H\grad f.
\end{array}\right.$$

\item If $M$ is a curve in $N(\alpha,\beta)$ with parallel mean curvature, then $M$ is $f$-biharmonic if and only if

$$ \left\{\begin{array}{l}

-f\tr B(\cdot,\nabla_{\cdot}\grad f)-f\tr\nabla_{\cdot}B(\cdot,\grad f)-n|\grad f|^2H-B(\grad f,\grad f)\\
=n^2f^2\alpha H-3nf^2\beta klH-nf^2\tr B(\cdot,A_H\cdot)+nf(\Delta f)H,\\ \\

+fRic_M(\grad f)+f\grad(\Delta f)+f\tr(A_{B(\cdot,\grad f)}(\cdot))-\dfrac12\grad(|\grad f|^2)\\
=2f(n-1)\alpha\grad f-3nfA_H\grad f.
\end{array}\right.$$

\end{enumerate}
\end{cor}

{\bf Proof:}
The proof is a direct consequence of Theorem \ref{thm1b} taking into account first that $M$ has parallel mean curvature so that the terms $\Delta^{\perp}H$, $\nabla^{\perp}_{\grad f}H$, $\grad|H|^2$ and $\tr(A_{\nabla^{\perp}_{\cdot}H\cdot})$ vanish. Moreover, we use the fact that
\begin{enumerate}
\item if $M$ is a hypersurface, then $m=0$ and so $jlH=0$, $kjH=0$ and $klH=-H$,
\item if $M$ is a complex surface then $k=0$ and $l=0$,
\item if $M$ is a Lagrangian surface, then $j=0$, $m=0$,
\item if $M$ is a curve, then $j=0$.
\end{enumerate}
\hfill $\square$
\begin{remark}
Note that from Theorem \ref{thm2b}, we can deduce a analogous corollary for hypersurfaces, curves and complex or Lagrangian submanifolds of complex space forms $\mathbb{M}^N_{\mathbb{C}}(c)$. Here again, the only difference is that $\alpha=\beta=c$ and $n<2N$ instead of $n<4$. We do not write down this corollary.
\end{remark}
\subsection{Bi-$f$-harmonic submanifolds of generalized Sasakian space forms}
\begin{thm}\label{thm3b}
Let $\psi: M^{p}\rightarrow \widetilde{M}(f_0,f_1,f_2)$ be a submanifold of a generalized Sasakian space form $\widetilde{M}(f_0,f_1,f_2)$, with second fundamental form $B$, shape operator $A$, mean curvature $H$ and a positive $C^{\infty}$-differentiable function $f$ on $M$. Then $M$ is $f$-biharmonic submanifold of $\widetilde{M}(f_0,f_1,f_2)$ if and only if the  following two equations are satisfied
\begin{enumerate}
\item
$$ \begin{array}{l}
nf^2\Delta^{\perp}H+nf^2\tr B(\cdot,A_H\cdot)-nf(\Delta f)H-3n\nabla^{\perp}_{\grad f}H\\
-f\tr B(\cdot,\nabla_{\cdot}\grad f)-f\tr\nabla_{\cdot}B(\cdot,\grad f)-n|\grad f|^2H-B(\grad f,\grad f)\\
=n^2f^2f_1H-nf^2f_2|\xi^{\top}|^2H-n^2f^2f_2\eta(H)\xi^{\perp}-3nf^2f_3NsH\\-(n-1)ff_2\eta(\grad f)\xi^{\perp}-3fNP\grad f,
\end{array}$$
\item $$\begin{array}{rcl}
&&\dfrac{n^2f^2}{2}\grad|H|^2+2n^2f^2\tr(A_{\nabla^{\perp}_{\cdot}H\cdot})+3nfA_H\grad f\\
&&+fRic_M(\grad f)+f\grad(\Delta f)+f\tr(A_{B(\cdot,\grad f)}(\cdot))-\dfrac12\grad(|\grad f|^2)\\
&=&-2n(n-1)ff_2\eta(H)\xi^{\top}-6nff_3PsH+(n-1)ff_1\grad f\\
&&-ff_2|\xi^{\top}|^2\grad f-(n-2)ff_2\eta(\grad f)\xi^{\top}-ff_3P^2\grad f.
\end{array}$$
\end{enumerate}

\end{thm}
{\bf Proof:} In order to prove this Theorem, we first recall that from the computation of Section \ref{sec3}, we have (see proof of Theorem \eqref{thm2})
\begin{equation}\label{eq3b0}
\tr\left(\overline{R}(\cdot,H)\cdot\right)^{\top}=f_2(n-1)\eta(H)\xi^{\top}+3f_3PsH
\end{equation}
and 
\begin{equation}\label{eq3b1}
\tr\left(\overline{R}(\cdot,H)\cdot\right)^{\perp}=-f_1nH+f_2\left(|\xi^{\top}|^2H+n\eta(H)\xi^{\perp}\right)+3f_3NsH.
\end{equation}
Moreover, we have this elementary lemma.
\begin{lemma}\label{lem3b2}
We have
$$\tr\left(\overline{R}(\cdot,\grad f)\cdot\right)^{\top}=-(n-1)f_1\grad f+f_2\left( |\xi^{\top}|^2 \grad f+(n-2)\eta(\grad f)\xi^{\top}\right)+f_3P^2\grad f$$
and 
$$\tr\left(\overline{R}(\cdot,\grad f)\cdot\right)^{\perp}=f_2(n-1)\eta(\grad f)\xi^{\perp}+3NP\grad f.$$
\end{lemma}
\hfill$\square$\\ \\
Now, combining  Equations \eqref{eq3b0}-\eqref{eq3b1} and Lemma \ref{lem3b2} together with Theorem \ref{thm1b}, we obtain\linebreak the conditions given in Theorem \ref{thm3b}, which concludes the proof.
\hfill$\square$\\ \\
Now, we give the proof of the Lemma \ref{lem3b2}.\\ \\
{\bf Proof of Lemma \ref{lem3b2}:} From the definition of the curvature tensor of $\widetilde{M}(f_1,f_2,f_3)$ we have
\begin{eqnarray*}
tr\left(\overline{R}(\cdot,\grad f)\cdot\right)&=&-(n-1)f_1\grad f +f_2\sum_{i=1}^n\Big(\eta(e_i)^2\grad f-\eta(\grad f)\eta(e_i)e_i\Big)\\
&&+f_2\sum_{i=1}^n\Big(<e_i,e_i>\eta(\grad f)\xi-<\grad f,e_i>\eta(e_i)\xi\Big)\\
&&+f_3\sum_{i=1}^n\Big(\Omega(e_i,\grad f)\phi e_i-\Omega)e_i,e_i)\phi\grad f+2\Omega(e_i,\grad f)e_i\\
&=&-(n-1)f_1\grad f +f_2\Big(|\xi^{\top}|^2\grad f-\eta(\grad f)\xi^{\top}+(n-1)\eta(\grad f)\xi\Big)\\
&&+3f_3\phi^2\grad f.
\end{eqnarray*}
Using the fact that $\phi^2\grad f=P^2\grad f+NP\grad f$ and identification of tangent and normal parts, we get the desired identities.
\hfill$\square$\\ \\
Here again, we finish this section with some particular cases. Namely, we have the \linebreak following corollary.

\begin{cor}\label{cor2b}
Let $\psi: M^{p}$ be a submanifold of a generalized Sasakian space form $\widetilde{M}(f_0,f_1,f_2)$ with parallel mean curvature. 
\begin{enumerate}
\item If $M$ is invariant then $M$ is bi-$f$-harmonic if and only if

$$ \left\{\begin{array}{l}
nf^2\tr B(\cdot,A_H\cdot)-nf(\Delta f)H\\
-f\tr B(\cdot,\nabla_{\cdot}\grad f)-f\tr\nabla_{\cdot}B(\cdot,\grad f)-n|\grad f|^2H-B(\grad f,\grad f)\\
=n^2f^2f_1H-nf^2f_2|\xi^{\top}|^2H-n^2f^2f_2\eta(H)\xi^{\perp}-(n-1)ff_2\eta(\grad f)\xi^{\perp},\\ \\

fRic_M(\grad f)+f\grad(\Delta f)+f\tr(A_{B(\cdot,\grad f)}(\cdot))-\dfrac12\grad(|\grad f|^2)\\
=-2n(n-1)ff_2\eta(H)\xi^{\top}-6nff_3PsH+(n-1)ff_1\grad f\\
-ff_2|\xi^{\top}|^2\grad f-(n-2)ff_2\eta(\grad f)\xi^{\top}-ff_3P^2\grad f- 2n^2f^2\tr(A_{\nabla^{\perp}_{\cdot}H\cdot}).
\end{array}\right.$$

\item  If $M$ is anti-invariant then $M$ is bi-$f$-harmonic if and only if

$$ \left\{\begin{array}{l}
nf^2\tr B(\cdot,A_H\cdot)-nf(\Delta f)H\\
-f\tr B(\cdot,\nabla_{\cdot}\grad f)-f\tr\nabla_{\cdot}B(\cdot,\grad f)-n|\grad f|^2H-B(\grad f,\grad f)\\
=n^2f^2f_1H-nf^2f_2|\xi^{\top}|^2H-n^2f^2f_2\eta(H)\xi^{\perp}-3nf^2f_3NsH\\-(n-1)ff_2\eta(\grad f)\xi^{\perp},\\ \\

fRic_M(\grad f)+f\grad(\Delta f)+f\tr(A_{B(\cdot,\grad f)}(\cdot))-\dfrac12\grad(|\grad f|^2)\\
=-2n(n-1)ff_2\eta(H)\xi^{\top}-2n^2f^2\tr(A_{\nabla^{\perp}_{\cdot}H\cdot})\\
-ff_2|\xi^{\top}|^2\grad f-(n-2)ff_2\eta(\grad f)\xi^{\top}.
\end{array}\right.$$

\item
If $\xi$ is normal to $M$ then $M$ is bi-$f$-harmonic if and only if

$$ \left\{\begin{array}{l}

-f\tr B(\cdot,\nabla_{\cdot}\grad f)-f\tr\nabla_{\cdot}B(\cdot,\grad f)-n|\grad f|^2H-B(\grad f,\grad f)\\
=n^2f^2f_1H-n^2f^2f_2\eta(H)\xi-3nf^2f_3NsH-(n-1)ff_2\eta(\grad f)\xi\\-nf^2\tr B(\cdot,A_H\cdot)+nf(\Delta f)H,\\ \\

fRic_M(\grad f)+f\grad(\Delta f)+f\tr(A_{B(\cdot,\grad f)}(\cdot))-\dfrac12\grad(|\grad f|^2)\\
=(n-1)ff_1\grad f-2n^2f^2\tr(A_{\nabla^{\perp}_{\cdot}H\cdot}).\\

\end{array}\right.$$

\item If $\xi$ is tangent to $M$ then $M$ is bi-$f$-harmonic if and only if

$$ \left\{\begin{array}{l}
nf^2\tr B(\cdot,A_H\cdot)-nf(\Delta f)H\\
-f\tr B(\cdot,\nabla_{\cdot}\grad f)-f\tr\nabla_{\cdot}B(\cdot,\grad f)-n|\grad f|^2H-B(\grad f,\grad f)\\
=n^2f^2f_1H-nf^2f_2H-3nf^2f_3NsH-3fNP\grad f,\\ \\

fRic_M(\grad f)+f\grad(\Delta f)+f\tr(A_{B(\cdot,\grad f)}(\cdot))-\dfrac12\grad(|\grad f|^2)\\
=-6nff_3PsH+(n-1)ff_1\grad f-2n^2f^2\tr(A_{\nabla^{\perp}_{\cdot}H\cdot})\\
-ff_2|\xi^{\top}|^2\grad f-(n-2)ff_2\eta(\grad f)\xi^{\top}-ff_3P^2\grad f.
\end{array}\right.$$

\item If $M$ is a hypersurface then $M$ is $f$-harmonic if and only if

$$ \left\{\begin{array}{l}
nf^2\tr B(\cdot,A_H\cdot)-nf(\Delta f)H\\
-f\tr B(\cdot,\nabla_{\cdot}\grad f)-f\tr\nabla_{\cdot}B(\cdot,\grad f)-n|\grad f|^2H-B(\grad f,\grad f)\\
=n^2f^2f_1H-nf^2f_2|\xi^{\top}|^2H-n^2f^2f_2\eta(H)\xi^{\perp}-(n-1)ff_2\eta(\grad f)\xi^{\perp}\\-3fNP\grad f,\\ \\

fRic_M(\grad f)+f\grad(\Delta f)+f\tr(A_{B(\cdot,\grad f)}(\cdot))-\dfrac12\grad(|\grad f|^2)\\
=-2n(n-1)ff_2\eta(H)\xi^{\top}+(n-1)ff_1\grad f-2n^2f^2\tr(A_{\nabla^{\perp}_{\cdot}H\cdot})\\
-ff_2|\xi^{\top}|^2\grad f-(n-2)ff_2\eta(\grad f)\xi^{\top}-ff_3P^2\grad f.
\end{array}\right.$$

\end{enumerate}
\end{cor}

{\bf Proof:} The proof is a direct consequence of Theorem \ref{thm3b} using the fact that the mean curvature is parallel and so the terms $\Delta^{\perp}H$, $\nabla^{\perp}_{\grad f}H$, $\grad|H|^2$ and $\tr(A_{\nabla^{\perp}_{\cdot}H\cdot})$ vanish. In addition, we use
\begin{enumerate}
\item if $M$ is invariant, then $P=0$,
\item if $M$ is anti-invariant, $N=0$,
\item if $\xi$ is normal, then $\eta(\grad f)=0$ and $M$ is anti-invariant which implies $P=0$,
\item if $\xi$ is tangent, then $\eta(H)=0$,
\item if $M$ is a hypersurface, then $sH=0$.
\end{enumerate}
\hfill $\square$

\section*{\textbf{Acknowledgements}}
Second author is supported by post doctoral scholarship of ``Harish Chandra Research Institute", Department of Atomic Energy, Government of India.


\begin{thebibliography}{00}
\bibitem{ABC} P. Alegre, D. E. Blair and A. Carriazo, {\it Generalized Sasakian space forms}, Israel J. Math., 141, 157 - 183, (2004).\\

\bibitem{AC} P. Alegre, and A. Carriazo, {\it Generalized Sasakian space forms and conformal change of the metric}, Results Math., 59(3), 485 - 493, (2011).\\

\bibitem{BFO} P. Baird, A. Fardoun and S. Ouakkas, {\it Conformal and semi-conformal biharmonic maps}, Ann. Glob. Anal. Geom., 34, 403 - 414, (2008).\\


\bibitem{BMO} A. Balmu\c{s}, S. Montaldo and C. Oniciuc, {\it On the biharmonicity of pseudo-umbilical and PNMC submanifolds in spheres and their type}, Ark. Mat., 51, 197 - 221, (2013).\\

\bibitem{Bla} D. E. Blair, {\it Riemannian Geometry of Contact and Symplectic Manifolds}, Birkh\"auser Boston, Progress in Mathematics, 203, (2002).\\

\bibitem{CMO} R. Caddeo, S. Montaldo, C. Oniciuc, {\it Biharmonic submanifolds in spheres.} Israel J. Math., 130, 109 - 123, (2002).\\

\bibitem{CMO1} R. Caddeo, S. Montaldo, and P. Piu, {\it On biharmonic maps}, Global differential geometry: the mathematical legacy of Alfred Gray, Contemp. Math., vol. 288, Amer. Math. Soc., Providence, RI, (2000).\\

\bibitem{Chen} B. Y. Chen, {\it Total Mean Curvature and Submanifolds of Finite Type}, Series in Pure Mathematics {\bf 1}, World Scientific Publishing Co., Singapore, (1984).\\

\bibitem{Ch2} B. Y. Chen, {\it Some open problems and conjectures on submanifolds of finite type}, Soochow J. Math., 17, \linebreak 169 - 188, (1991).\\

\bibitem{Ch3} B. Y. Chen, {\it Recent developments of biharmonic conjecture and modified biharmonic conjectures}, Pure and Applied Differential Geometry, Proceedings of the conference PADGE 2012, Shaker Verlag, Aachen, 81 - 90, (2013).\\

\bibitem{Der} A. Derdzinski, {\it Exemples de m\'etriques de K\"ahler et d'Einstein auto-duales sur le plan complexe}, G\'eom\'etrie riemannienne en dimension 4, (S\'eminaire Arthur Besse 1978/79), Cedic/Fernand Nathan, Paris, 334 - 346, (1981).\\ 

\bibitem{ES} J. Eells and J. H. Sampson, {\it Harmonic mappings of Riemannian manifolds}, Amer. J. Math., 86, 109 - 160, (1964).\\

\bibitem{EL} J. Eells and L. Lemaire, {\it Selected topics in harmonic maps}, CBMS, 50, Amer. Math. Soc, (1983).\\

\bibitem{FLMO} D. Fetcu, E. Loubeau, S. Montaldo and C. Oniciuc, {\it Biharmonic submanifolds of $\CC P^n$}, Math. Z., 266, 505 - 531, (2010).\\

\bibitem{FO} D. Fetcu, C. Oniciuc, {\it Explicit formulas for biharmonic submanifolds in Sasakian space forms},  Pacific J. Math., 240 (1), 85 - 107, (2009).\\

\bibitem{FOR} D. Fetcu, C. Oniciuc and H. Rosenberg, {\it Biharmonic submanifolds with parallel mean curvature in $\SSS^n\times\RR$}, J. Geom. Anal. (in press).\\

\bibitem{JI} J. Inoguchi, {\it Submanifolds with harmonic mean curvature vector field in contact 3-manifolds}, Colloq. Math., 100, 163 - 179, (2004).\\

\bibitem{Ji} G. Y. Jiang, {\it 2-harmonic maps and their first and second variational formulas}, Chinese Ann. Math. Ser., A7 (4), 389 - 402, (1986).\\

\bibitem{Lot} A. Lotta, {\it Slant submanifolds in contact geometry}, Bull. Math. Soc. Roumanie, 39, 183 - 198, (1996).\\

\bibitem{LO} T. Liang and Y. -L. Ou, {\it Biharmonic hypersurfaces in a conformally flat space}, Results Math. 64, 91 - 104, (2013). \\

\bibitem{LM} E. Loubeau and S. Montaldo, {\it Biminimal immersions}, Proc. Edinb. Math. Soc., 51, 421 - 437, (2008).\\

\bibitem{LU} Wei-Jun Lu, {\it On  f-Biharmonic  maps  between  Riemannian  manifolds}, arXiv:1305.5478, preprint, (2013).\\

\bibitem{SH} S. Maeta and H. Urakawa, {\it Biharmonic Lagrangian submanifolds in Kaehler manifolds}, Glasgow Math. J., 55, 465 - 480, (2013).\\

\bibitem{MC} S. Montaldo and C. Oniciuc, {\it A short survey on biharmonic maps between Riemannian manifolds}, Rev. Un. Mat. Argentina, 47, 1 - 22 (2007).\\

\bibitem{OND}   S. Ouakkas, R. Nasri, and M. Djaa, {\it On the f-harmonic and f-biharmonic maps}, JP J. Geom. Topol. 10 (1), 11-27, (2010).\\

\bibitem{Ols} Z. Olszak, {\it On the existence of generalized complex space forms}, Israel J. Math., 65, no. 2, 214 - 218, (1989).\\

\bibitem{OL} Y. -L. Ou and L. Tang, {\it The generalized ChenÕs conjecture on biharmonic submanifolds is false}, Michigan Math. J., 61, 531 - 542, (2012).\\

\bibitem{YO} Y. -L. Ou, {\it Biharmonic hypersurfaces in Riemannian manifolds}, Pacific J. Math., Vol. 248, No. 1, 217 - 232, (2010).\\

\bibitem{YO1} Y. -L. Ou, {\it Some recent progress of biharmonic submanifolds}, arXiv:1511.09103.\\

\bibitem{YO2} Y. -L. Ou, {\it On f-biharmonic maps and f-biharmonic submanifolds}, Pacific J. Math., 271, 461-477, (2014).\\

\bibitem{YO3} Y. -L. Ou, {\it f-Biharmonic maps and f-biharmonic submanifolds II},  arXiv:1605.00128.\\

\bibitem{Rot} J. Roth, {\it A note on biharmonic submanifolds of product spaces}, J. Geom., 104, 375 - 381, (2013). \\

\bibitem{RotAu} J. Roth and A. Upadhyay, {\it Biharmonic submanifolds of generalized space forms}, arXiv:1602.06131.\\ 

\bibitem{TV} F. Tricerri and L. Vanhecke, {\it Curvature tensors on almost Hermitian manifolds},  Trans. Amer. Math. Soc., 267, no. 2, 365 - 397, (1981).\\

\bibitem{HU} H. Urakawa, {\it Sasaki manifolds, K\"ahler cone manifolds and biharmonic submanifolds}, Illinois J. Math., Volume 58, Number 2, 521 - 535, (2014).\\

\bibitem{YK} K. Yano and M. Kon, {\it Structures on manifolds}, Series in Pure Mathematics, 3. World Scientific Publishing Co., Singapore, (1984).\\

\end{thebibliography}
\end{document}